\documentclass[11pt,leqno]{amsart}
\usepackage{amsthm,amsfonts,amssymb,amsmath,oldgerm,mathabx}
\numberwithin{equation}{section}
\usepackage{fullpage}
\usepackage{color}

\usepackage[pdftex]{graphicx} 
\usepackage{hyperref}
\usepackage{stmaryrd}
\usepackage{comment}

\usepackage{caption}
\usepackage{subcaption}

\newcommand{\SortNoop}[1]{}

\renewcommand\d{\partial}
\DeclareMathOperator\eD{e}
\DeclareMathOperator\iD{i}

\def\eps{\varepsilon }
\DeclareMathOperator{\Id}{Id}
\DeclareMathOperator{\I}{I}
\DeclareMathOperator{\b0}{\mathbf{0}}
\newcommand{\transp}[1]{{#1}^{{\sf T}}}
\newcommand{\cond}[1]{{\bf (D#1)}}

\DeclareMathOperator{\Span}{Span}

\DeclareMathOperator{\dD}{d}
\DeclareMathOperator{\Div}{\transp{\nabla}}

\DeclareMathOperator{\beD}{\mathbf{e}}
\DeclareMathOperator{\supp}{supp}

\newcommand\br{\begin{remark}}
\newcommand\er{\end{remark}}
\newcommand\bp{\begin{pmatrix}}
\newcommand\ep{\end{pmatrix}}
\newcommand{\be}{\begin{equation}}
\newcommand{\ee}{\end{equation}}
\newcommand{\bes}{\begin{equation*}}
\newcommand{\ees}{\end{equation*}}
\newcommand\ba{\begin{equation}\begin{aligned}}
\newcommand\ea{\end{aligned}\end{equation}}
\newcommand\bas{\begin{equation*}\begin{aligned}}
\newcommand\eas{\end{aligned}\end{equation*}}

\newcommand\nn{\nonumber}
\newcommand{\beg}{\begin{example}}
\newcommand{\eeg}{\end{example}}
\newcommand{\bpr}{\begin{proposition}}
\newcommand{\epr}{\end{proposition}}
\newcommand{\bt}{\begin{theorem}}
\newcommand{\et}{\end{theorem}}
\newcommand{\bc}{\begin{corollary}}
\newcommand{\ec}{\end{corollary}}
\newcommand{\bl}{\begin{lemma}}
\newcommand{\el}{\end{lemma}}
\newcommand{\bd}{\begin{definition}}
\newcommand{\ed}{\end{definition}}
\newcommand{\brs}{\begin{remarks}}
\newcommand{\ers}{\end{remarks}}

\newtheorem{theorem}{Theorem}[section]
\newtheorem{proposition}[theorem]{Proposition}
\newtheorem{corollary}[theorem]{Corollary}
\newtheorem{lemma}[theorem]{Lemma}

\theoremstyle{remark}
\newtheorem{remark}[theorem]{Remark}
\theoremstyle{definition}
\newtheorem{definition}[theorem]{Definition}

\newtheorem{example}[theorem]{Example}

\newcommand\R{\mathbf R}
\newcommand\C{\mathbf C}
\newcommand{\N}{\mathbf N}
\newcommand{\Z}{\mathbf Z}

\newcommand\bA{{\mathbf A}}

\newcommand\bD{{\mathbf D}}

\newcommand\bF{{\mathbf F}}
\newcommand\bG{{\mathbf G}}

\newcommand\bK{{\mathbf K}}
\newcommand\bL{{\mathbf L}}
\newcommand\bM{{\mathbf M}}

\newcommand\bP{{\mathbf P}}
\newcommand\bQ{{\mathbf Q}}

\newcommand\bU{{\mathbf U}}
\newcommand\bV{{\mathbf V}}
\newcommand\bW{{\mathbf W}}
\newcommand\bX{{\mathbf X}}
\newcommand\bY{{\mathbf Y}}

\newcommand\bDelta{\boldsymbol{\Delta}}

\newcommand\bLambda{\boldsymbol{\Lambda}}

\newcommand\bPhi{\boldsymbol{\Phi}}
\newcommand\bPsi{\boldsymbol{\Psi}}
\newcommand\bOm{\boldsymbol{\Omega}}

\newcommand\bfc{{\mathbf c}}

\newcommand\bff{{\mathbf f}}
\newcommand\bfg{{\mathbf g}}

\newcommand\bfp{{\mathbf p}}
\newcommand\bfq{{\mathbf q}}

\newcommand\bfx{{\mathbf x}}

\newcommand\bfxi{{\boldsymbol \xi}}

\newcommand\bfphi{\boldsymbol{\phi}}
\newcommand\bfvarphi{\boldsymbol{\varphi}}

\newcommand\ubK{{\underline \bK}}

\newcommand\ubM{{\underline \bM}}

\newcommand\ubU{{\underline \bU}}

\newcommand\ubc{{\underline \bfc}}

\newcommand\tbP{\widetilde{\bP}}
\newcommand\tbQ{\widetilde{\bQ}}

\newcommand\tbU{\widetilde{\bU}}

\newcommand\tbW{\widetilde{\bW}}

\newcommand\tbY{\widetilde{\bY}}

\newcommand\tchi{\widetilde{\chi}}

\newcommand\tbq{\widetilde{\bfq}}

\newcommand\uM{{\underline M}}

\newcommand\uOm{{\underline \Omega}}

\newcommand\uc{{\underline c}}

\newcommand\cA{{\mathcal A}}
\newcommand\cB{{\mathcal B}}
\newcommand\cC{{\mathcal C}}
\newcommand\cD{{\mathcal D}}
\newcommand\cE{{\mathcal E}}
\newcommand\cF{{\mathcal F}}

\newcommand\cJ{{\mathcal J}}
\newcommand\cK{{\mathcal K}}

\newcommand\cM{{\mathcal M}}
\newcommand\cN{{\mathcal N}}
\newcommand\cO{{\mathcal O}}

\newcommand\cR{{\mathcal R}}
\newcommand\cS{{\mathcal S}}

\newcommand\cU{{\mathcal U}}

\newcommand\cW{{\mathcal W}}

\newcommand\bcK{\boldsymbol{\cK}}

\newcommand\bcM{\boldsymbol{\cM}}

\newcommand\bcU{\boldsymbol{\cU}}

\newcommand\bcW{\boldsymbol{\cW}}

\usepackage{enumitem}

\title[2D periodic waves in systems with conservation laws]{Nonlinear stability of two-dimensional periodic waves in parabolic systems with conservation laws}

\author{L.~Miguel Rodrigues}
\address{
Univ Rennes, CNRS, IRMAR - UMR 6625, F-35000 Rennes, France}
\email{{\tt luis-miguel.rodrigues@univ-rennes1.fr}}
\thanks{Research of L.M.R. was supported by the ANR Project HEAD ANR-24-CE40-3260, the Institut Universitaire de France, and by France 2030 through the programme Centre Henri Lebesgue ANR-11-LABX-0020-01.}

\author{Aric Wheeler}
\address{Duke University, 27708 Durham, North Carolina}
\email{{\tt aaw74@math.duke.edu}}
\thanks{Research of A.W. was supported by NSF grants DMS-1700279 and DMS-2038056.}

\begin{document}

\begin{abstract}
We develop a stability theory for two-dimensional periodic traveling waves of general parabolic systems, possibly including conservation laws. In particular, we identify a diffusive spectral stability assumption and prove that it implies nonlinear stability for variously-(non)localized perturbations, including critically nonlocalized perturbations. Thus we extend the stability parts of \cite{JNRZ-conservation} to two-dimensional patterns and of \cite{MR} to systems with conservation laws. In doing so we need to bypass two kinds of low spectral regularity, explicitly conic-like singularities due to multidimensionality and Jordan-block like singularities due to conservation laws.

\vspace{0.5em}

{\small \paragraph {\bf Keywords:} periodic traveling-wave solutions; parabolic systems; asymptotic stability.
}

\vspace{0.5em}

{\small \paragraph {\bf AMS Subject Classifications:} 35B35, 35K40, 35C07, 35B40, 35B10, 37L15.
}
\end{abstract}

\date{\today}
\maketitle

\tableofcontents

\section{Introduction}\label{s:introduction}

We continue here the general programme, initiated in \cite{MR}, aiming at a complete stability theory for genuinely multi-dimensional periodic traveling waves of parabolic systems, extending to the multidimensional context the comprehensive theory available for plane periodic waves \cite{JNRZ-conservation}.

We carry out most of our present analysis in the two-dimensional context but discuss higher-dimensional extensions in Appendix~\ref{s:3D}. Namely we consider essentially
\be\label{eq:rdc}
\bcW_t=\Delta \bcW+\Div\bG(\bcW)+\bff(\bcW)\,,
\ee
for the $\R^n$-valued unknown $\bcW$, $\bcW(t,\bfx)\in\R^n$ (with $n\in\N^\star$), $t$ denoting the time variable and $\bfx\in\R^2$ the spatial variables. In \eqref{eq:rdc}, we identify elements of $\R^n$ with column vectors, the operator $\Div$ acts row-wise and $\Delta$ is the scalar Laplacian. The flux function $\bG:\R^n\to \cM_{2,n}(\R)$ and the reaction rate $\bff:\R^n\to\R^n$ are assumed to be smooth. Finally, we assume that the first $r$ components of $\bff$ vanish identically, so that the first $r$ equations of \eqref{eq:rdc} are conservation laws.

Our goal is to study the dynamics of \eqref{eq:rdc} near one of its (uniformly) traveling solutions, that is near a solution of the form 
\[
\bcW(t,\bfx)=\bcU(\bfx-t\,\bfc)\,,
\]
where $\bfc\in\R^2$ is the wavespeed vector and $\bcU$ is the wave profile of the wave under consideration. We restrict to two-dimensional periodic traveling waves, meaning that there exists a fixed basis $(\bX_1,\bX_2)$ of $\R^2$ such that the wave profile $\bcU$ is left invariant by translations in both the $\bX_1$ and $\bX_2$ directions. Alternatively, letting $(\bK_1,\bK_2)$ denote the dual basis of $(\bX_1,\bX_2)$, one may express such traveling waves in the form
\begin{equation}\label{eq:alternativeform}
		\bcW(t,\bfx)=\bU(\transp{\bK}(\bfx-t\,\bfc))=\bU(\transp{\bK}\bfx+t\,\bOm),
\end{equation}
where $\bK=\bp \bK_1&\bK_2 \ep\in\cM_{2,2}(\R)$ is the matrix of wave vectors, $\bOm=-\transp{\bK}\bfc$ is the time frequency vector and the scaled wave profile $\bU$ is left invariant by translations by the canonical basis of $\R^2$. Throughout we identify waves by their mathematical structure but, to ease comparison with the rest of the literature, we stress that they commonly receive more vivid names, planar periodic waves being often designed as 
stripes or rolls whereas two-dimensional periodic waves include cases designed as square, hexagon, herringbone, \emph{etc.} patterns, depending on the geometric shape of $(\bX_1,\bX_2)$. For a more thorough discussion on the existence of waves we simply refer to the introductions of \cite{MR,BdR-R}.

\subsection*{Bloch symbols}

The fact that waves are not isolated but form a smoothly parameterized family plays a deep role in the organization of the dynamics. When we pick a single specific wave we shall materialize it by underlining the various features of the specific wave, $\ubU$, $\uc$, $\uOm$, \emph{etc.} To study its stability it is expedient to move to an adapted co-moving frame
\begin{align}\label{eq:co-moving}
\bcW(t,\bfx)\,=\,\bW\left(t,\transp{\ubK}\,(\bfx-t\,\ubc)\right)
\end{align}
so that \eqref{eq:rdc} is turned into
\be\label{eq:rd}
\bW_t=\transp{(\ubK\nabla)}(\ubK\nabla)\bW +\transp{(\ubK\nabla)} \bG(\bW) + (\transp{\ubK}\ubc \cdot \nabla) \bW +\bff(\bW).
\ee
admitting $\ubU$ as a stationary $\Z^2$-periodic solution. Linearizing \eqref{eq:rd} about $\ubU$ yields the periodic-coefficient equation $(\d_t-L)\bV=0$ with $L$ given by
\[
L\bV:= \transp{(\ubK\nabla)}(\ubK\nabla) \bV + \transp{(\ubK\nabla)} \dD\bG(\ubU)(\bV) + (\transp{\ubK}\ubc \cdot \nabla) \bV + \dD\bff(\ubU)(\bV)\,.
\]
The operator $L$ is seen as acting on $L^2(\R^2;\C^n)$ with domain $H^2(\R^2;\C^n)$. The spectral stability of $\ubU$ is then expressed in terms of the spectrum of $L$. 

As classical for periodic-coefficient operators the spectrum of $L$ is more conveniently handled by introducing a suitable integrable transform --- the Floquet-Bloch transform --- so as to study $L$ through its Bloch symbols $L_\bfxi$, parametrized by Floquet exponents $\bfxi\in[-\pi,\pi]^2$. Each $L_\bfxi$ acts on $L^2([0,1]^2;\C^n)\cong L^2(\R^2/\Z^2;\C^n)$ with domain $H^2_{\rm per}([0,1]^2;\C^n)\cong H^2(\R^2/\Z^2;\C^n)$ through 
\[
L_{\bfxi}\bV=
\transp{(\ubK(\nabla+\iD\bfxi))}\ubK(\nabla+\iD\bfxi)\bV +\transp{(\ubK(\nabla+\iD\bfxi))} \dD\bG(\ubU)(\bV) + (\transp{\ubK}\ubc \cdot (\nabla+\iD\bfxi))\bV + \dD\bff(\ubU)(\bV)\,,
\]
thus has compact resolvents hence discrete spectrum, reduced to eigenvalues of finite multiplicity. In particular, the following spectral decomposition holds
\[
\sigma(L)\,=\,\bigcup_{\bfxi\in[-\pi,\pi]^2} \sigma(L_\bfxi)\,.
\]
We recall the associated minimal background in Subsection~\ref{s:Bloch} and refer to \cite[Appendix~A]{MR} for a more detailed treatment.

\subsection*{Jordan block structure}

The presence of conservation laws in \eqref{eq:rdc} deeply affects the nature of the spectrum of $L_{\b0}$ compared to the situation dealt with in \cite{MR}, essentially as the situation in \cite{JNRZ-conservation} differs from the one in \cite{JNRZ-RD1,JNRZ-RD2}. The generalized kernel of $L_{\b0}$ is directly related to the properties of periodic traveling waves near $\ubU$ that share the same lattice of periods. In the absence of conservation laws, as in systems of \cite{MR}, generically the set of nearby co-periodic waves is reduced to $\ubU$ and its spatial translates and, accordingly, the generalized kernel of $L_{\b0}$ is reduced to its kernel and spanned by spatial derivatives of $\ubU$, thus here two-dimensional and spanned by $\d_1\ubU$, $\d_2\ubU$. In turn, the fact that the first $r$ equations of \eqref{eq:rdc} are conservation laws directly implies that the kernel of the adjoint of $L_{\b0}$ contains the functions constant equal to $\beD_j$, $j=1,\cdots,r$ (where $(\beD_1,\cdots,\beD_n)$ denotes the canonical basis of $\R^n$). Since those are independent and orthogonal to $\d_1\ubU$ and $\d_2\ubU$, the generalized kernel of $L_{\b0}$ must be of dimension at least $r+2$.

The following assumptions thus encode generic criticality (in the stable case).
\begin{enumerate}
  \item[\cond1] The spectrum of $L_{\b0}$ intersects $\iD\R$ only at $\lambda=0$, $\lambda=0$ being an eigenvalue of $L_{\b0}$ of algebraic multiplicity $r+2$. 
  \item[\cond2] For any $\bfxi\neq\b0$, $\sigma(L_{\bfxi})\subset\{\,\lambda\,;\,\Re(\lambda)<0\,\}$.
\end{enumerate}

Under Condition~\cond1, one derives that, near $\ubU$, periodic traveling waves of wavevector matrix $\bK$ take the form $(\bU,\bfc)=(\bU^{\bK,\bM}(\cdot+\bfvarphi),\bfc^{\bK,\bM})$, where $\bfvarphi\in\R^2$ and $\bM\in\R^r$ with $M_j=\int_{[0,1]^2}\beD_j\cdot \bU$. In particular, $(\d_1\ubU,\d_2\ubU,\d_{M_1}\ubU,\cdots,\d_{M_r}\ubU)$ span the generalized kernel of $L_{\b0}$ with 
\[
L_{\b0}(\d_{M_j}\ubU)=-(\transp{\ubK}\d_{M_j}\ubc \cdot \nabla)\ubU\,.
\]
Therefore, unless $\dD_{\bM}\ubc$ vanishes, the generalized kernel of $L_{\b0}$ contains nontrivial Jordan blocks (as in \cite{JNRZ-conservation}). Before stating our main results we still need to explain how those Jordan blocks affect the dynamics.

\subsection*{Space-modulated stability}

In \cite{JNRZ-conservation} the main dynamical effect was the necessity of proving stability in a space-modulated sense. Translated into our current multi-D framework this consists in measuring in the classical definition of stability the proximity to $\ubU$ of both initial data and solutions at later times by 
\[
\inf_{\bPhi\textrm{ invertible}}\|\bW\circ\bPhi-\ubU\|_{\bX} + \|\nabla(\bPhi-\Id)\|_{\bY}
\]
where $(\bX,\bY)$ are some functional spaces (possibly different for initial data and the solution). This is crucial when $\bPhi-\Id$ is not decaying with time but $\nabla(\bPhi-\Id)$ is. For more informal insights we refer to \cite{R,R_Roscoff}. For comparison, note that naive stability requires control on $\|\bW-\ubU\|_{\bX}$ whereas orbital stability requires control on 
\[
\inf_{\bPhi\textrm{ uniform translation}}\|\bW\circ\bPhi-\ubU\|_{\bX}\,.
\]
Note that it is possible to choose $\bY$ as a space of curl-free vector-fields, on which it is natural to choose $\|\,\cdot\,\|_\bY$ as $\|\bV\|_{\bY}:=\|\Div(\bV)\|_{\tbY}$ for some functional space $\tbY$. 

The notion encodes that one should introduce a \emph{local} phase shift to capture the main effect of perturbations and that replacing this local phase shift with its derivative tames the growth induced by Jordan blocks. It is intuitively clear that when $\dD_{\bM}\ubc$ is nonzero, local perturbations in shape will result in modifications in local speeds which in turn will accumulate over time into local modification of phase shifts. Formal geometrical optics expansions in the spirit of \cite{Noble-Rodrigues} give more quantitative educated guesses.

The upshot of such a formal analysis is that the large-time dynamics is expected to be described at leading order by modulation in parameters
\begin{align*}
\bcW(t,\bfx)&\sim \bU^{(\bcK,\bcM)(t,\bfx)}(\bPsi(t,\bfx))\,,&
\bcK&=\nabla\bPsi\,,
\end{align*}
with $(\bcK,\bcM)$ remaining close to $(\ubK,\ubM)$ (in the stable case) and solving 
\begin{align}\label{eq:W}
\d_t\bp\bcK\\\bcM\ep
&=\bp\nabla (\bOm^{(\bcK,\bcM)})\\ \Div\bF^{(\bcK,\bcM)} \ep+\bLambda[\nabla](\bcK,\bcM)\,,&\bcK\quad\textrm{curl-free}
\end{align}
where $\bOm^{(\bcK,\bcM)}=-\transp{\bcK}\,\bfc^{(\bcK,\bcM)}$, $\bF^{(\bcK,\bcM)}\in \cM_{2,r}$, $\bF^{(\bcK,\bcM)}_{j,\ell}:=\int_{[0,1]^2} \bG_{j,\ell}\circ \bU^{(\bcK,\bcM)}$, and $\bLambda[\nabla]$ is a (less explicit) second-order operator compatible with the propagation of the curl-free condition.

\subsection*{Diffusive spectral stability}

We shall make only a loose use of the foregoing expected modulation behavior. The main insight we want to extract is how the spectral singularity resulting from Jordan blocks should be blown up. For this purpose we may safely replace the curl-free wavevector variable $\bcK=\nabla \bPsi$ with $(-\Delta)^{1/2}\bPsi$ and drop the curl-free constraint, since $-(-\Delta)^{-1/2}\Div$ and $\nabla (-\Delta)^{-1/2}$ may be used as $L^2$ isometric inverses. 

To give more details about the connection let us first introduce, for $\bfxi$ sufficiently small, a basis $(\bQ_1^\bfxi,\cdots,\bQ_{r+2}^\bfxi)$ of the sum of generalized eigenspaces of $L_\bfxi$ associated with its $r+2$ eigenvalues near $0$, starting from $(\d_1\ubU,\d_2\ubU,\d_{M_1}\ubU,\cdots,\d_{M_r}\ubU)$ at $\bfxi=\b0$, and a dual basis $(\tbQ_1^\bfxi,\cdots,\tbQ_{r+2}^\bfxi)$ playing the same role for the adjoint of $L_\bfxi$, starting from\footnote{We identify here constant functions with their values.} $(*,*,\beD_1,\cdots,\beD_r)$ at $\bfxi=\b0$. In particular for $\bfxi$ small the spectrum of $L_\bfxi$ near $0$ coincides with the spectrum of the matrix
\[
\bD_\bfxi:=(\langle \tbQ_j^\bfxi;L_\bfxi \bQ_\ell^\bfxi\rangle_{L^2([0,1]^2)})_{1\leq j,\ell\leq r+2}\,.
\]
Note that
\[
\bD_{\b0}=\bp \b0_{2,2}& -\transp{\ubK}\dD_\bM\ubc\\\b0_{r,2}&\b0_{r,r}\ep\,.
\]

By design, in the chosen basis the first two coordinates play the role of (the Fourier transform of) local phases. To mimick the change $\bPsi\mapsto (-\Delta)^{1/2}\bPsi$, we introduce for $\bfxi\neq\b0$ small
\[
\bDelta_\bfxi:=\bp \|\bfxi\|\,\I_2&\b0_{2,r}\\\b0_{r,2}&\I_r\ep\,\bD_\bfxi\,\bp \|\bfxi\|\,\I_2&\b0_{2,r}\\\b0_{r,2}&\I_r\ep^{-1}
\]
extended as $\bDelta_{\b0}=\b0_{r+2,r+2}$. We prove in Subsection~\ref{s:side-band} that despite the apparent singularity the map $\bfxi\mapsto \bDelta_\bfxi$ is indeed smooth in polar coordinates. 

The following assumption encodes that System~\eqref{eq:W} is diffusive near $(\ubK,\ubM)$ (at least at low-frequency).

\begin{enumerate}
  \item[\cond3] There exist positive $\theta$, $C$ and $\xi_0$ such that when $\|\bfxi\|\leq \xi_0$, the matrix $\bDelta_\bfxi$ defined above (using Condition~\cond1) satisfies for any $t\geq0$
\[
\|\eD^{t\,\bDelta_\bfxi}\|\leq C\,\eD^{-t\,\theta\,\|\bfxi\|^2}\,.
\]
\end{enumerate}

The control afforded by Conditions~\cond1-\cond2-\cond3 is weaker than the ones provided by diffusive stability in both \cite{JNRZ-conservation} and \cite{MR} with many respects. In all cases diffusive stability implies that there exists a $\theta'>0$ such that for any $\bfxi$
\[
\sigma(L_\bfxi)\subset\{\ \lambda\,;\,\Re(\lambda)\leq-\theta'\|\bfxi\|^2\ \}\,.
\]
In \cite{MR} the absence of conservation laws thus of Jordan blocks at the critical part of the spectrum allowed the latter spectral assumption to be reinforced into: for some positive $C'$, $\theta'$, for any $\bfxi$,
\[
\|e^{tL_\bfxi}\|_{L^2([0,1]^2)\to L^2([0,1]^2)}\leq C'\,\eD^{-t\theta'\|\bfxi\|^2},
\]
where $\eD^{tL_\bfxi}$ denotes the semigroup generated by the Bloch symbol $L_\bfxi$. A similar bound holds for the one-dimensional conservationless case \cite{JNRZ-RD1} but it comes as a consequence of the spectral bound alone, whereas in the multi-dimensional case the bound contains extra information; see the detailed discussion in \cite[Appendix~A]{MR} and \cite{BdR-R}. In contrast, here, as in \cite{JNRZ-conservation}, diffusive stability only yields: for some positive $C'$, $\theta'$, for any $\bfxi$,
\[
\|e^{tL_\bfxi}\|_{L^2([0,1]^2)\to L^2([0,1]^2)}\leq C'\,\frac{\eD^{-t\theta'\|\bfxi\|^2}}{\|\bfxi\|}\,.
\]

In the one-dimensional case one may benefit from the fact that one-parameter families of operators have much better spectral regularity properties. In particular, in \cite{JNRZ-conservation}, beyond the desingularization of phase/wavenumber type turning $\bD_\bfxi$ into $\bDelta_\bfxi$, one may also enforce that $\bfxi\mapsto \bDelta_\bfxi$ is smooth\footnote{This requires us to use at the modulation level $\nabla \bPsi$ and not $(-\Delta)^{1/2}\bPsi$, and at spectral level to replace $\|\bfxi\|$ with $\bfxi$ in the definition of $\bDelta_\bfxi$, which is only possible in dimension $1$.} and not only smooth in polar coordinates, and that it is smoothly diagonalized. The distinction here is similar to the one for first-order hyperbolic systems: whereas in one dimension the strict hyperbolicity may be thought as generic, there is a much larger variety in higher dimensions. In particular the acoustic system --- that is, the system form of the wave equation --- share almost no common feature with a system of uncoupled transport equations in dimension at least two. The formulation of Condition~\cond3 is precisely designed to not dive into this complexity and encompass all the reasonable cases without distinction.

Part of this multiD system complexity is also present in \cite{MR}. However the modulation system there, replacing \eqref{eq:W}, has only two equations. The most refined analysis in \cite{MR} uses the special structure of systems of two equations. Note for instance that hyperbolic constant-coefficient systems of two first-order equations are either strictly hyperbolic or consist of two uncoupled scalar equations. In contrast we have here $r+2$ equations to deal with and the situation with three equations in 2D is expected to be close to generic. To give a hint in this direction, we point out that with three equations in dimension two one may already design a system that is hyperbolic in each direction but not hyperbolic ; see Petrowski's example in \cite{Benzoni-Serre}.

The main technical question our current analysis answers positively is whether with both difficulties jointly occurring one may still conclude to nonlinear stability.

\subsection*{Main results}

Before stating our main results the last thing remaining to discuss is which class of perturbations we are considering. The level of localization of the initial perturbation is directly connected to the rate of time decay one may expect. On general grounds, with this is mind, there are two natural choices. On one hand one may pick the maximal level of localization that may be converted into time decay, thus obtaining the maximal time decay (for general\footnote{Commonly one may beat this decay rate on a subset of positive codimension, being defined by some exceptional vanishing. For instance mean-free solutions to the heat equation decay faster than the heat decay rate.} data). On the other hand one may pick the localization resulting in the minimal time decay rate compatible with a proof of nonlinear stability by the Duhamel formula, encoding the variation of constant strategy. The latter situation is said to be critical and it is characterized by the fact that in the Duhamel formula the linear and nonlinear part decay at the same rate, whereas in subcritical situations the nonlinear part decays faster.

Here the multidimensionality actually helps since the decay mechanism is spreading (by diffusion) and is thus stronger in higher dimensions. In the one-dimensional case with conservation laws of \cite{JNRZ-conservation} there is no choice: maximal decay and critical decay are at the same level. It corresponds to taking an initial phase shift $\bfphi_0$ such that $\nabla \bfphi_0$ is integrable. The presence of a Jordan block at the spectral origin prevents further localization to be of any use. The critical conservationless one-dimensional case is analyzed in \cite{JNRZ-RD1}. Yet for the same class of equations the maximal decay case was analyzed decades before in the pioneering \cite{Schneider}. The approach of Schneider sets a renormalization procedure, thus relies heavily on the self-similar character of the Burgers equation (playing the role of \eqref{eq:W}) and is therefore by design restricted to the conservationless one-dimensional case. It was subsequently extended to the critical decay case in \cite{SSSU}. Proofs of \cite{JNRZ-RD1,JNRZ-conservation} rather build from the revisiting of the maximal decay rate in \cite{JZ}. As an interesting recent piece of work we point out \cite{dR} revisiting the crictical conservationless one-dimensional case, that combines arguments similar to those of \cite{JNRZ-RD1} with special properties of the Burgers equation so as to replace the integrability assumption on $\nabla \bfphi_0$ with a boundedness assumption on $\bfphi_0$, both being, in dimensional one, essentially\footnote{We use the homogeneity of seminorms under the action of homotethies as a measure of how much regularity/localization the seminorm controls. In dimension $d$, $s$ derivatives in $L^p$ is associated with the index $s-d/p$, the lower it is the more localization and the less regularity we control.} at the same level of localization but the latter is slightly more general than the former.

In the two-dimensional conservationless framework treated in \cite{MR} the maximal decay threshold is at the integrable $\bfphi_0$ level and the critical threshold is at the integrable $\Delta\bfphi_0$ level. In the present framework the presence of a Jordan block at the spectral origin reduces the range of possibilities from integrable $\nabla\bfphi_0$ to integrable $\Delta\bfphi_0$.

\bt[Critical decay]\label{theo:critmain}
Assume \cond1-\cond2-\cond3. There exist $\eps_0>0$ and $C>0$ such that if for some sublinear\footnote{By this, we mean as in \cite{MR} that $\bfphi_0$ may differ from $\Delta^{-1}(\Delta \bfphi_0)$ by a constant function but not by a non-constant affine function. This is for instance the case if one enforces $\nabla\bfphi_0\in\Span\left(\bigcup_{1\leq p<\infty}L^p(\R^2;\cM_2(\R))\right)$.} $\bfphi_0$
\[
E_0:=\|\bW_{0}(\cdot-\bfphi_{0})-\ubU \|_{(H^{2} \cap W^{2,4})(\R^2;\R^n)}+\|\Delta \bfphi_{0}\|_{(L^{1} \cap W^{1,4})(\R^2;\R^2)}\,\leq\,\eps_0
\]
then, there exist a unique global solution to \eqref{eq:rd} with initial datum $\bW_{0}$ and a phase shift $\bfphi$ with $\bfphi(0,\cdot)= \bfphi_{0}$ such that, for any $t\geq0$,
\begin{align*}
\|\bW(t,\cdot-\bfphi(t,\cdot))-\ubU\|_{W^{2,4}(\R^2;\R^n)}
+\|\nabla\bfphi(t,\cdot)\|_{W^{2,4}(\R^2;\cM_2(\R))}
+\|\d_t\bfphi(t,\cdot)\|_{W^{2,4}(\R^2;\R^2)}
&\leq \frac{C\,E_{0}}{(1+t)^{\frac{1}{4}}}\,.
\end{align*}
Furthermore, with constants independent of $(\bW_0,\bfphi_0)$ and no further restriction on $E_0$, 
\begin{enumerate}
\item for any $t\geq0$, 
\begin{align*}
\|\bW(t,\cdot-\bfphi(t,\cdot))-\ubU\|_{L^\infty(\R^2;\R^n)}
+\|\nabla\bfphi(t,\cdot)\|_{L^\infty(\R^2;\cM_2(\R))}
+\|\d_t\bfphi(t,\cdot)\|_{L^\infty(\R^2;\R^2)}
&\leq C\,E_0\,\frac{\ln(2+t)}{(1+t)^{\frac{1}{2}}}\,;
\end{align*}
\item for any $2<p_{0}<q_{0}<\infty$, there exists a constant $C_{p_{0},q_{0}}>0$, such that for any $p \in [p_{0},q_{0}]$, and any $t\geq0$
\[
\|\bW(t,\cdot-\bfphi(t,\cdot))-\ubU\|_{L^p(\R^2;\R^n)}
+\|\nabla\bfphi(t,\cdot)\|_{L^p(\R^2;\cM_2(\R))}
+\|\d_t\bfphi(t,\cdot)\|_{L^p(\R^2;\R^2)}
\,\leq \frac{C_{p_{0},q_{0}}\,E_{0}}{(1+t)^{\frac{1}{2}-\frac{1}{p}}}\,.
\]
\end{enumerate}
\et

\bt[Maximal decay]\label{theo:submain}
Assume \cond1-\cond2-\cond3. There exist $\eps_0>0$ and $C>0$ such that if for some $\bfphi_0$
\[
\cE_{0}:=\|\bW_{0}(\cdot-\bfphi_{0})-\ubU \|_{(L^{1}\cap H^2\cap W^{2,4})(\R^2;\R^n)}+\|\nabla \bfphi_{0}\|_{(L^{1} \cap H^2\cap W^{2,4})(\R^2;\cM_2(\R))}\,\leq\,\eps_0
\]
then, there exist a unique global solution to \eqref{eq:rd} with initial datum $\bW_{0}$ and a phase shift $\bfphi$ with $\bfphi(0,\cdot)= \bfphi_{0}$ such that, for any $t\geq0$, for any $2\leq p\leq 4$
\begin{align*}
\|\bW(t,\cdot-\bfphi(t,\cdot))-\ubU\|_{W^{2,p}(\R^2;\R^n)}
+\|\nabla\bfphi(t,\cdot)\|_{W^{2,p}(\R^2;\cM_2(\R))}
+\|\d_t\bfphi(t,\cdot)\|_{W^{2,p}(\R^2;\R^2)}
&\leq \frac{C\,\cE_{0}}{(1+t)^{1-\frac{1}{p}}}\,.
\end{align*}
Furthermore, with no further restriction
\begin{enumerate}
\item for any $2 \leq p \leq \infty$ and any $t\geq0$
\[
\|\bW(t,\cdot-\bfphi(t,\cdot))-\ubU\|_{L^p(\R^2;\R^n)}
+\|\nabla\bfphi(t,\cdot)\|_{L^p(\R^2;\cM_2(\R))}
+\|\d_t\bfphi(t,\cdot)\|_{L^p(\R^2;\R^2)}
\,\leq \frac{C\,\cE_{0}}{(1+t)^{1-\frac{1}{p}}}\,,
\]
\item there exists a constant $\bfphi_{\infty}$ depending only on $\bfphi_{0}$ such that for any $t \geq 0$
\[
\|\bW(t,\cdot)-\ubU \|_{L^{p}(\R^2;\R^n)} + \| \bfphi(t,\cdot)- \bfphi_{\infty} \|_{L^{p}(\R^2;\R^2)} \leq \frac{C \cE_{0}}{(1+t)^{\frac{1}{2}-\frac{1}{p}}}  \,, \qquad 2  \leq p \leq +\infty\,.
\]
\end{enumerate}
\et

Having System~\eqref{eq:W} in mind, the foregoing decay rates should be compared with the decay rates of the two-dimensional heat equation in the variable $\nabla \bfphi$. Note however that we only state and prove such rates in $L^p$-based spaces with $2\leq p\leq \infty$. As may also be guessed from \eqref{eq:W} this is not a purely technical restriction. The first-order part of the modulation system if it contains parts of wave-type induces dispersive effects that deteriorate decay rates in $L^p$, $p<2$. We refer the reader to the detailed analysis of the two-dimensional barotropic compressible Navier-Stokes system in \cite{Hoff_Zumbrun-NS_compressible_pres_de_zero,Rodrigues-compressible} for a sharp related discussion.

Those dispersive effects are precisely the spatial counterparts to the low regularity on the spectral side, due to multidimensionality. In particular our main challenge is to check that the analysis of \cite{JNRZ-conservation} may be extended to the present case without suffering from $L^p$, $p<2$, deterioration. To stress that the outcome is not completely obvious, we point out that \cite{Hoff_Zumbrun-NS_compressible_pres_de_zero,Rodrigues-compressible} deals with the maximal decay case whereas the critical decay case for small solutions of the two-dimensional barotropic compressible Navier-Stokes system is widely open; see the discussion in \cite{Rodrigues-these}. 

Even when dealing with the critical decay case we benefit from the fact that non critical linear estimates do hold in order to deal with nonlinear contributions. Note that this does not hold for \cite{JNRZ-conservation} but the dispersive obstruction is absent in dimension one. The benefit is similar (but smaller) than the one one uses to close convective nonlinear stability with no decay at all, as for instance in \cite{GR}.

We postpone to Subsection~\ref{s:critmain} a thorough technical discussion of the scheme of the proof of Theorem~\ref{theo:critmain}. Yet we would like to stress that its proof differs significantly from proofs in both \cite{JNRZ-conservation} and \cite{MR}. When examining the contribution of nonlinear terms in a Duhamel formula only terms with a good conservative structure may be treated as in \cite{MR}. Unfortunately the introduction of the phase does bring new non-conservative terms. To deal with those one must use and prove that higher derivatives of the phase decay faster. This is expected from the linear analysis but a proof at the nonlinear level requires the propagation of those higher decay rates in a critical decay regime. The foregoing is inherently difficult. A similar difficulty occurs in \cite{JNRZ-conservation} but it is bypassed there by resorting to estimates on the linear propagator that are precluded here by the possible dispersive effects. Instead, we show that for a suitable Bloch multiplier operator $\cJ$, $\cJ\bV$ also decays faster at the linear level and prove the nonlinear propagation of those decay rates jointly with those of the higher derivatives of the phase. This requires a (soft) multilinear analysis adapted to $\cJ$.

\subsection*{Perspectives}

In \cite{JNRZ-RD2,JNRZ-conservation,MR} is also provided a detailed description of the large-time asymptotic behavior including the validation of a system such as \eqref{eq:rd}. We expect that some form of validation could also be obtained here but a rough computation, in the critical case, only shows that refinement in modulation (in all parameters) yield an improvement of the decay rates of the remainders in $L^p$ from $(1+t)^{-(1/2-1/p)}$ to $(1+t)^{-1/2}$. It is far from obvious that the later rate is optimal and we have decided to save some technicalities and leave this for further investigation. As in the current stability analysis the technical issue is that in general, unlike in \cite{MR}, one faces dispersion in some of the components but not in all of them so that one may have to pay the price of worst rates in $L^p$, $p<2$ without benefiting from extra decay in $L^p$, $p>2$.

In order to apply the present results and those in \cite{MR} it is important to provide stability diagrams in all relevant bifurcation scenarios, multiplying the type of analysis carried out in \cite{BdR-R}. The analysis there untangles for conservationless equations the situation where plane waves become unstable by a transverse instability at nonzero frequencies. The stable two-dimensional patterns then appear as secondary instablities and the arising modulation system is indeed fully dispersive. We stress that even the case when small-amplitude two-dimensional patterns arise from primary instability of constant solutions needs further investigation. Indeed whereas the corresponding existence part has been thoroughly developed for a long time --- as far ago as \cite{Sattinger_bifurcation} ---, the accompanying stability is currently restricted at best\footnote{For patterns with extra symmetry, such as hexagonal pattern, the class of studied perturbations is further restricted by the extra symmetries.} to co-periodic perturbations, corresponding to Floquet $\bfxi=\b0$. We see the development of a complete bifurcation analysis incorporating diffusive spectral stability as the primary remaining target to reach.

We expect that as in \cite{BdR-R} the elucidation of bifurcations will also shed some light on the structure of modulation systems near bifurcation and provide sufficient elements to clarify the analysis of asymptotic behavior. To be somewhat more concrete, in scenarios\footnote{We expect the latter for instance when $r=1$.} when the modulation system would exhibit a clear separation between a dispersive part and a scalar-type part one could prove that at leading order only the dispersiveless part matter in large-time. This would echo the fact that solutions to the barotropic compressible Navier-Stokes equations become asymptotically incompressible in large-time \cite{Hoff_Zumbrun-NS_compressible_pres_de_zero,Rodrigues-compressible}.

\subsection*{Outline}

After the present introduction the rest of the paper is organized as follows. In the next section we gather some preliminary results (elements of Bloch-wave spectral analysis, geometric structure of profile equations, desingularization by phase modulation). Then in the two sections after that, we provide first linear estimates then nonlinear estimates for the critical case, thus in the end proving Theorem~\ref{theo:critmain}. In the final section before the appendix, we prove Theorem~\ref{theo:submain}. Finally, in the appendix we discuss extensions to higher dimensions, covering the critical three-dimensional case.

\medskip

\noindent\emph{Acknowledgment.} A.W. thanks the University of Rennes for its hospitality when the project was initiated. Both authors are grateful to Kevin Zumbrun for his continued interest in this project.

\section{Preliminaries}\label{s:prelim}

\subsection{Basic facts about the Floquet-Bloch analysis}\label{s:Bloch}

We gather here a minimal set of useful facts about the Floquet-Bloch transform. As a preliminary we make explicit our convention for the Fourier transform. When $g\in L^1(\R^2)$, its Fourier transform is defined by
\[
\cF(g)(\bfxi)=\widehat{g}(\bfxi)\,:=\,\frac{1}{(2 \pi)^2} \int_{\R^2} \eD^{-\iD \bfxi \cdot \bfx} g(\bfx)\,\dD\bfx\,.
\]
Then when $g\in \cS(\R^2)$ its Floquet-Bloch transform is defined by
\[
\cB(g)(\bfxi,\bfx)
=\widecheck{g}(\bfxi,\bfx)\,:=\,
\sum_{\bfp \in \Z^2} \eD^{2\iD\pi \bfp \cdot \bfx} 
\widehat{g}(\bfxi + 2\pi \bfp)
\,=\,
\frac{1}{(2\pi)^{2}} \sum_{\bfq \in \Z^2} \eD^{-\iD \bfxi \cdot(\bfx + \bfq)} g(\bfx + \bfq)\,,
\]
the equivalence of both formula being a special case of the Poisson summation formula. Note that $\cB(g)(\bfxi,\cdot)$ is $(\beD_1,\beD_2)$-periodic for any $[-\pi,\pi]^2$.

\bpr
Let $g\in \cS(\R^2)$. The following properties hold.
\begin{enumerate}
\item \emph{Floquet-Bloch inversion.} For any $\bfx\in\R^2$
\[
g(\bfx) = \int_{[-\pi,\pi]^{2}} \eD^{\iD \bfxi \cdot \bfx}\,\widecheck{g}(\bfxi,\bfx) \dD\bfxi\,.
\]
\item \emph{Differentiation.} For any $(\bfxi,\bfx)\in [-\pi,\pi]^2\times\R^2$
\[
\widecheck{(\nabla_\bfx g)}(\bfxi,\bfx)\,=\,(\nabla_\bfx+\iD\bfxi)(\widecheck{g})(\bfxi,\bfx)\,.
\]
\item \emph{Periodic multiplication.} When $h$ is $(\beD_1,\beD_2)$-periodic, for any $(\bfxi,\bfx)\in [-\pi,\pi]^2\times\R^2$
\[
\widecheck{(g\,h)}(\bfxi,\bfx)\,=\,h(\bfx)\,\widecheck{g}(\bfxi,\bfx)\,.
\]
\item \emph{Low-frequency functions.} If $\supp\widehat{g}\subset[-\pi,\pi]^2$, for any $(\bfxi,\bfx)\in [-\pi,\pi]^2\times\R^2$
\[
\widecheck{g}(\bfxi,\bfx)\,=\,\widehat{g}(\bfxi)\,.
\]
\item \emph{Parseval identity.} 
\[
\|g\|_{L^2(\R^2)}\,=\,(2\pi)\|\widecheck{g}\|_{L^{2}([-\pi,\pi]^{2};L^{2}([0,1]^2))}\,.
\]
\item \emph{Hausdorff-Young inequalities.} For any $2\leq p\leq \infty$ and $p'$ its Lebesgue conjugate,\\
that is, $\frac{1}{p}+\frac{1}{p'}=1$ 
\begin{align*}
\|g\|_{L^{p}(\R^2)}&\leq (2\pi)^{\frac{2}{p}} \|\widecheck{g}\|_{L^{p'}([-\pi,\pi]^{2};L^{p}([0,1]^2))}\,,\\
 \|\widecheck{g}\|_{L^p([-\pi,\pi]^2;L^{p'}([0,1]^2))}&\leq (2\pi)^{-\frac{2}{p'}} 
 \|g\|_{L^{p'}(\R^2)}\,.
\end{align*}
\end{enumerate}
\epr

Properties $(2)$ and $(3)$ imply the announced $\cB(L\bfg)(\bfxi,\bfx)=L_{\bfxi}(\cB(\bfg)(\bfxi,\cdot))(\bfx)$. Properties $(3)$ and $(4)$ explain how by taking a Floquet-Bloch transform and averaging one may extract the evolution of parameters, such as phase shifts, out of linear slow modulations. Properties $(5)$ and $(6)$ are key to extend the Floquet-Bloch transform by duality and continuity.

Let us denote with an index ${\cdot}_{\rm per}$ the closure of (restrictions of) $\cC^\infty$ functions that are $(\beD_1,\beD_2)$-periodic.

\bpr
\begin{enumerate}
\item $(2\pi)\cB$ is a one-to-one isometry from $L^2(\R^2)$ to $L^{2}([-\pi,\pi]^{2};L^{2}([0,1]^2))$.
\item For any $s\in\R$, $\cB$ is one-to-one between $H^s(\R^2)$ and $L^2([-\pi,\pi]^2;H^s_{\rm per}([0,1]^2))$.
\item For any $s\in\R_+$, $2\leq p\leq \infty$, $\cB$ is continuous from $L^{p'}([-\pi,\pi]^2;W^{s,p}_{\rm per}([0,1]^2))$ to $W^{s,p}(\R^2)$ and $\cB^{-1}$ is continuous from $W^{s,p'}(\R^2)$ to $L^{p}([-\pi,\pi]^2;W^{s,p'}_{\rm per}([0,1]^2))$.
\end{enumerate}
\epr

Finally we introduce the Bloch multiplier operator $\cJ$ defined by
\[
\cB(\cJ g)(\bfxi,\bfx):=\bfxi\,\cB(g)(\bfxi,\bfx)\,.
\]
The convention is that it acts component-wise on vector-valued functions. Note that $\cJ$ commutes with any operator defined by Bloch symbols, such as $L$, or in other words with any operator diagonalized by the Floquet-Bloch transform in the sense that it acts Floquet exponent by Floquet exponent. We shall use $\cJ$ to track that in the dynamical situation we consider a higher cancellation at low-Floquet exponents is associated with a faster time decay.

The following lemma proves a Leibniz' rule for multilinear expressions of sufficiently low-Floquet functions, particularly convenient to deal with nonlinear terms.

\bl\label{l:Leibniz}
Let $m\in\N$, $m\geq 2$. If the Floquet-Bloch transforms of $g_1$, $\cdots$, $g_m$ are supported where $\|\bfxi\|\leq \pi/m$, then for any $\alpha\in\N^2$
\[
\cJ^\alpha(g_1\cdots g_m)\,=\,
\sum_{\substack{(\alpha_{(1)},\cdots,\alpha_{(m)})\in(\N^2)^m\\ \alpha_{(1)}+\cdots+\alpha_{(m)}=\alpha}}
\bp \alpha\\ \alpha_{(1)},\cdots,\alpha_{(m)}\ep\,\cJ^{\alpha_{(1)}}(g_1)\cdots \cJ^{\alpha_{(m)}}(g_m)\,.
\]
\el

\begin{proof}
This follows by combining that from any $(\zeta_{(1)},\cdots,\zeta_{(m)})\in(\R^2)^m$
\[
(\zeta_{(1)}+\cdots+\zeta_{(m)})^\alpha\,=\,
\sum_{\substack{(\alpha_{(1)},\cdots,\alpha_{(m)})\in(\N^2)^m\\ \alpha_{(1)}+\cdots+\alpha_{(m)}=\alpha}}
\bp \alpha\\ \alpha_{(1)},\cdots,\alpha_{(m)}\ep\,\zeta_{(1)}^{\alpha_{(1)}}\cdots \zeta_{(m)}^{\alpha_{(m)}}
\]
and that from the low-Floquet assumption stems
\[
\cB(g_1\cdots g_m)(\bfxi,\bfx)
\,=\,\left(\cB(g_1)(\cdot,\bfx)\star\cdots\star\cB(g_m)(\cdot,\bfx)\right)(\bfxi)\,.
\]
\end{proof}

\subsection{Profile variations}\label{s:profile}

We now untangle the structure of nearby waves under assumption \cond{1}. In doing so we mostly follow \cite[Appendix~B1]{MR}.

Let $\ubM\in\R^r$ be defined by $\uM_j=\int_{[0,1]^2}\beD_j\cdot \ubU$ and observe that the profile equation for waves of wavematrix $\bK$ and speed $\bfc$ is 
\be\label{stand_eq}
0=\transp{(\bK\nabla)}(\bK\nabla)\bU +\transp{(\bK\nabla)} \bG(\bU) + (\transp{\bK}\bfc \cdot \nabla) \bU +\bff(\bU).
\ee

\bpr\label{p:structure}
Assume~\cond{1}. Then there exist $\eps_0>0$ and a smooth map 
\[
B((\ubK,\ubM),\eps_0)\to H^2_{\rm per}([0,1]^2;\R^n)\times\R^2\,,\quad
(\bK,\bM)\mapsto (\bU^{\bK,\bM}(\cdot),\bfc(\bK,\bM))
\]
such that, for any $(\bK,\bM)\in B(\ubK,\eps_0)$, $(\bK,\bU^{\bK,\bM}(\cdot),\bfc(\bK,\bM))$ solves \eqref{stand_eq} with 
\begin{align}\label{mean-norm}
M_j&=\int_{[0,1]^2}\beD_j\cdot \bU^{\bK,\bM}\,,& j=1,\cdots,r\,,
\end{align}
and for any $(\bU,\bfc)\in H^2_{\rm per}([0,1]^2;\R^n)\times\R^2$ such that $(\bU,\bfc)$ solves \eqref{stand_eq}-\eqref{mean-norm} and 
\begin{align*}
\|\bfc-\ubc\|&\leq \eps_0\,,&
\inf_{\bfvarphi_0\in\R^2}\|\bU-\ubU(\cdot+\bfvarphi_0)\|_{H^2_{\rm per}([0,1]^2;\R^n)}&\leq \eps_0\,,
\end{align*}
one has $\bfc=\bfc(\bK,\bM)$ and there exists $\bfvarphi\in\R^2$ such that $\bU=\bU^{\bK,\bM}(\cdot+\bfvarphi)$. Moreover the map $(\bK,\bM)\mapsto \bU^{\bK,\bM}$ is valued in $H^\infty_{\rm per}([0,1]^2;\R^n)$ and, for any $s\in\N$, there exists $0<\eps_0'\leq\eps_0$ such that it is smooth as a map from $B((\ubK,\ubM),\eps_0')$ to $H^s_{\rm per}([0,1]^2;\R^n)$. 
\epr

\begin{proof}
The proof follows the Lyapunov-Schmidt reduction. Let $\Pi_{\b0}$ denote the spectral projector associated with the eigenvalue $0$ of $L_{\b0}$ and complete first $(*,*,\beD_1,\cdots,\beD_r)$ into a basis $(\tbQ_1^{\b0},\cdots,\tbQ_{r+2}^{\b0})$ of the kernel of $\Pi_{\b0}^*$ such that 
\begin{align*}
\langle \tbQ_j^{\b0}; \d_\ell\ubU\rangle_{L^2([0,1]^2;\R^n)}&\,=\,\delta_{j,\ell}\,,&
1\leq j,\ell\leq2\,,
\end{align*}
then $(\d_1\ubU,\d_2\ubU,*,\cdots,*)$ into a basis $(\bQ_1^{\b0},\cdots,\bQ_{r+2}^{\b0})$ of the kernel of $\Pi_{\b0}$ dual to $(\tbQ_1^{\b0},\cdots,\tbQ_{r+2}^{\b0})$. To enforce uniqueness we shall reduce to seeking an unknown $\bV$ lying in the kernel of $\Pi_{\b0}$, that is, orthogonal to $(\tbQ_1^{\b0},\cdots,\tbQ_{r+2}^{\b0})$. 

We first show that we can factor out translational invariance and enforce orthogonality to $(\tbQ_1^{\b0},\tbQ_2^{\b0})$. To do so, we may apply the Implicit Function Theorem to the map 
\[
H^2([0,1]^2;\R^n)\times\R^2\to\R^2\,,\qquad(\bU,\bfvarphi)\longrightarrow
(\langle \tbQ_j^{\b0};\bU(\cdot-\bfvarphi)\rangle_{L^2([0,1]^2;\R^n)})_{j=1,2}\,.
\]
Indeed the map is $\cC^1$ and, at $(\ubU,\b0)$, its differential map with respect to $\bfvarphi$ is $-\I$. By using translational invariance, this implies that there exist $\eps>0$ and $C>0$ such that if $(\bU,\bfvarphi_0)$ is such that
\[
\|\bU-\ubU(\cdot+\bfvarphi_0)\|_{H^2_{\rm per}([0,1]^2;\R^n)}\leq \eps
\]
then there exists $\bfvarphi$ such that $\tbU=\bU(\cdot-\bfvarphi)$ satifies
\begin{align}\nonumber
\|\tbU-\ubU\|_{H^2_{\rm per}([0,1]^2;\R^n)}&\leq C\,\|\bU-\ubU(\cdot+\bfvarphi_0)\|_{H^2_{\rm per}([0,1]^2;\R^n)}\,,&\\\label{orth}
\langle \tbQ_j^{\b0};(\tbU-\ubU)\rangle_{L^2([0,1]^2;\R^n)}&=0\,,&1\leq j\leq 2\,.
\end{align}

It is thus sufficient to prove genuine uniqueness under the assumption that $\bU-\ubU$ is small and satisfies \eqref{orth}. Let us denote by $L_{\b0}^\dagger$ the inverse of $L_{\b0}$ restricted to the range of $(\I-\Pi_{\b0})$. With the extra constraint \eqref{orth}, Equations~\eqref{stand_eq}-\eqref{mean-norm} are equivalent to
\[
\bU-\ubU\,=\,\bp \bQ_3^{\b0}&\cdots&\bQ_{r+2}^{\b0}\ep \bM+\bV
\]
with
\begin{align*}
\bV&\,=\,-L_{\b0}^\dagger[(\I-\Pi_{\b0})\,\cR]\\
\bfc&\,=\,\transp{(\bK^{-1}\ubK)}\ubc
-\transp{(\bK^{-1})}
\,\bp \langle\tbQ_1^{\b0};\cdot\rangle_{L^2([0,1]^2;\R^n)}\\
\langle\tbQ_2^{\b0};\cdot\rangle_{L^2([0,1]^2;\R^n)}\ep\,\cR
\end{align*}
where
\begin{align*}
\cR&=\left(\transp{(\bK\nabla)}(\bK\nabla)-\transp{(\ubK\nabla)}(\ubK\nabla)\right)\bU
+((\transp{\bK}\bfc-\transp{\ubK}\ubc)\cdot \nabla)(\bU-\ubU)\\
&\quad+\transp{(\bK\nabla)}\bG(\bU)
-\transp{(\ubK\nabla)}\left(\bG(\ubU)+\dD \! \bG(\ubU)(\bU-\ubU)\right)\\
&\quad+\bff(\bU)-\bff(\ubU)-\dD \! \bff(\ubU)(\bU-\ubU))\,.
\end{align*}
Note that in the foregoing reformulation we are using that from the conservative structure of the original system stems that $\Pi_{\b0}\cR=0$ is equivalent to $\langle \tbQ_j^{\b0};\cR\rangle_{L^2([0,1]^2;\R^n)}=0$, $1\leq j\leq 2$. The proof is then achieved by another application of the Implicit Function Theorem solving for $(\bV,\bfc)$ such that $(\I-\Pi_{\b0})\,\bV=0$.
\end{proof}

Let us stress that the above construction also yields $(\d_{M_1}\ubU,\cdots,\d_{M_r}\ubU)=(\bQ_3^{\b0},\cdots,\bQ_{r+2}^{\b0})$.

\subsection{Side-band structure}\label{s:side-band}

We turn to discuss the desingularization of the Jordan block structure at $0$ when $\bfxi=\b0$. By standard spectral perturbation theory, for which we refer to \cite[Appendix~A2]{MR}, we obtain, for $\bfxi$ sufficiently small, 
\begin{itemize}
\item a spectral projector $\Pi_\bfxi$, smoothly dependent on $\bfxi$, associated with the spectrum of $L_\bfxi$ near the origin
\item a basis $(\bQ_1^\bfxi,\cdots,\bQ_{r+2}^\bfxi)$ of the range of $\Pi_\bfxi$ smoothly dependent on $\bfxi$ and extending the above $(\bQ_1^{\b0},\cdots,\bQ_{r+2}^{\b0})$ 
\item a dual basis $(\tbQ_1^\bfxi,\cdots,\tbQ_{r+2}^\bfxi)$ of the range of $\Pi_\bfxi^*$ smoothly dependent on $\bfxi$ and extending the above $(\tbQ_1^{\b0},\cdots,\tbQ_{r+2}^{\b0})$.
\end{itemize}
In particular with 
\[
\bD_\bfxi:=(\langle \tbQ_j^\bfxi;L_\bfxi \bQ_\ell^\bfxi\rangle_{L^2})_{1\leq j,\ell\leq r+2}
\]
there holds for $\bfxi$ sufficiently small
\[
\eD^{t\,L_\bfxi}\Pi_\bfxi
\,=\,
\bp \bQ_1^{\bfxi}&\cdots&\bQ_{r+2}^{\bfxi}\ep
\eD^{t\,\bD_\bfxi}
\bp \langle\tbQ_1^{\bfxi};\cdot\rangle_{L^2}\\
\vdots\\
\langle\tbQ_{r+2}^{\bfxi};\cdot\rangle_{L^2}\ep\,.
\]

To prepare further analysis of $\bD_\bfxi$, we expand
\[
L_\bfxi \bV\ =\ L_{\b0}\bV\ +\ \transp{(L^{(1)}\bV)}\iD(\ubK\bfxi)
\ -\ \|\ubK\bfxi\|^2\, \bV
\]
where
\[
L^{(1)}\bV \ :=\ 2 \ubK\nabla \bV + \dD \! \bG(\ubU)(\bV) + \ubc \transp{\bV}.
\]

In order to blow up the Jordan block, we introduce for $\bfxi\neq\b0$ small
\[
\bDelta_\bfxi:=\bp \|\bfxi\|\,\I_2&\b0_{2,r}\\\b0_{r,2}&\I_r\ep\,\bD_\bfxi\,\bp \|\bfxi\|\,\I_2&\b0_{2,r}\\\b0_{r,2}&\I_r\ep^{-1}
\]
extended as $\bDelta_{\b0}=\b0_{r+2,r+2}$. Note that by setting
\begin{align*}
\bp \bfq_1^{\bfxi}&\cdots&\bfq_{r+2}^{\bfxi}\ep
&:=\,\bp \bQ_1^{\bfxi}&\cdots&\bQ_{r+2}^{\bfxi}\ep
\,\bp \I_2&\b0_{2,r}\\\b0_{r,2}&\|\bfxi\|\,\I_r\ep\\
\bp \tbq_1^{\bfxi}&\cdots&\tbq_{r+2}^{\bfxi}\ep
&:=\,\bp \tbQ_1^{\bfxi}&\cdots&\tbQ_{r+2}^{\bfxi}\ep
\,\bp \|\bfxi\|\,\I_2&\b0_{2,r}\\\b0_{r,2}&\I_r\ep
\end{align*}
one obtains for $\bfxi$ sufficiently small
\[
\eD^{t\,L_\bfxi}\Pi_\bfxi
\,=\,
\frac{1}{\|\bfxi\|}\bp \bfq_1^{\bfxi}&\cdots&\bfq_{r+2}^{\bfxi}\ep
\eD^{t\,\bDelta_\bfxi}
\bp \langle\tbq_1^{\bfxi};\cdot\rangle_{L^2}\\
\vdots\\
\langle\tbq_{r+2}^{\bfxi};\cdot\rangle_{L^2}\ep\,.
\]

The following proposition is the main result of the present subsection.

\bpr\label{smooth-delta}
Under \cond{1}, the map $\bfxi\mapsto \bDelta_\bfxi$ is indeed smooth in polar coordinates. 
\epr

\begin{proof}
Note that since $\bQ_1^{\b0}$, $\bQ_2^{\b0}$ lie in the kernel of $L_{\b0}$ and $\tbQ_3^{\b0}$, $\cdots$, $\tbQ_{r+2}^{\b0}$ lie in the kernel of $L_{\b0}^*$, we already know that
\begin{align*}
\bD_{\b0}&\,=\,\bp \b0_{2,2}&*\\\b0_{r,2}&\b0_{r,r}\ep\,,&
\bD_{\bfxi}&\,=\,\bp \cO(\|\bfxi\|)&\cO(1)\\\cO(\|\bfxi\|)&\cO(\|\bfxi\|)\ep\,,
\end{align*}
so that there only remains to show that for any $\bfxi$
\begin{align*}
(\dD_{\bfxi}\bD_{\b0}(\bfxi))_{j,\ell}&=0\,,&
(\bD_{\bfxi})_{j,\ell}&=\cO(\|\bfxi\|^2)\,,&
3\leq j\leq r+2\,,\ 1\leq \ell\leq2\,.
\end{align*}
For the same structural reason, for such a $(j,\ell)$,
\[
(\dD_{\bfxi}\bD_{\b0}(\bfxi))_{j,\ell}
\,=\,
\langle \tbQ_j^{\b0};\transp{(L^{(1)}\bQ_\ell^{\b0})}\iD(\ubK\bfxi)\rangle_{L^2}
\,=\,
\langle \beD_{j-2};\transp{(L^{(1)}\d_\ell\ubU)}\iD(\ubK\bfxi)\rangle_{L^2}\,.
\]
At last, from the conservative structure of the system and the explicit formula of $L^{(1)}$ stem the claimed cancellation. 
\end{proof}

Though we have decided to prove Proposition~\ref{smooth-delta} in a more elementary way, we'd like to stress that it follows in a more robust way by differentiating along the family of periodic traveling waves. See \cite{BGNR,KR,Audiard-Rodrigues} for similar results in more elaborate situations.

From the definition of $\bDelta_{\bfxi}$ one also derives readily the following proposition. 

\bpr\label{p:Bloch-bounds}
Under assumptions \cond{1}-\cond{2}-\cond{3}, there exist positive $\xi_0$, $\theta'$ and $C'$ such that the following hold.
\begin{enumerate}
\item For any $\|\bfxi\|\geq \xi_0$
\[
\|e^{tL_\bfxi}\|_{L^2\to L^2}\leq C'\,\eD^{-t\theta'}\,.
\]
\item For any $\|\bfxi\|\leq \xi_0$
\begin{align*}
\|e^{tL_\bfxi}\,(\I-\Pi_\bfxi)\|_{L^2\to L^2}&\leq K\,\eD^{-t\theta'}\,,\\
\|e^{t\bD_\bfxi}\|_{\C^{r+2}\to \C^{r+2}}&\leq \frac{C'}{\|\bfxi\|}\,\eD^{-t\theta'\|\bfxi\|^2}\,,\\
\|e^{tL_\bfxi}\,\Pi_\bfxi\|_{L^2\to L^2}&\leq \frac{C'}{\|\bfxi\|}\,\eD^{-t\theta'\|\bfxi\|^2}\,.
\end{align*}
\end{enumerate}
\epr

\begin{proof}
The last two estimates follow from the definition of $\bDelta_{\bfxi}$ whereas the two first ones are derived from the fact that $L_\bfxi$ when $\|\bfxi\|\geq\xi_0$ and the restriction of $L_\bfxi$ to the kernel of $\Pi_\bfxi$ when $\|\bfxi\|\geq\xi_0$ generate analytic semigroups and possess uniform spectral gaps. 
\end{proof}

When analyzing the phase contributions the following lemma will also turn out to be useful. It plays a role similar to \cite[Lemma~3.1]{MR} and a similar normalization is also hidden in the proof of \cite[Proposition~1.7]{JNRZ-conservation}.

\bl\label{lem:blocksplit}
Under assumption~\cond{1}, one may also enforce that for $3\leq j\leq r+2$, $1\leq \ell\leq2$,
\begin{align*}
\langle \tbQ_j^\bfxi;\bQ_\ell^{\b0}\rangle_{L^2}&=\cO(\|\bfxi\|^2)\,,&
\langle \tbQ_j^{\b0};\bQ_\ell^\bfxi\rangle_{L^2}&=\cO(\|\bfxi\|^2)\,.&
\end{align*}
\el

\begin{proof}
Both terms are smooth and vanish at $\bfxi=\b0$ so that we only need to check that their linear expansions vanish identically. Since for such $(j,\ell)$, $\langle \tbQ_j^\bfxi;\bQ_\ell^{\bfxi}\rangle_{L^2}=0$, the linear expansion of the first term is equal to the opposite of the linear expansion of the second term thus they vanish simultaneously.

To conclude the proof, we only need to prove that by replacing $((\bQ^{\bfxi}_1,\cdots,\bQ^{\bfxi}_{r+2}))_{\bfxi}$, $((\tbQ^{\bfxi}_1,\cdots,\tbQ^{\bfxi}_{r+2}))_{\bfxi}$ with some $((\bP^{\bfxi}_1,\cdots,\bP^{\bfxi}_{r+2}))_{\bfxi}$, $((\tbP^{\bfxi}_1,\cdots,\tbP^{\bfxi}_{r+2}))_{\bfxi}$ satisfying the same spectral conditions one may also achieve the extra normalization condition: for any $\bfxi$, and $(j,\ell)$ as above 
\[
\langle \tbP_\ell^{\b0};
\dD_{\bfxi} \bP^{\b0}_{j}(\bfxi)\rangle_{L^2}
\,=\,0\,.
\]
This may be achieved, for $\bfxi$ sufficiently small, through
\begin{align*}
\bP_\ell^{\bfxi}&:=
\bQ_\ell^{\bfxi}
-\sum_{j=3}^{r+2} \langle \tbQ_j^{\b0};\dD_{\bfxi} \bQ^{\b0}_\ell(\bfxi)\rangle_{L^2}\,\bQ_j^{\bfxi}\,,&
1\leq\ell\leq2\,,\\
\bP_\ell^{\bfxi}&:=
\bQ_\ell^{\bfxi}\,,&
3\leq \ell\leq r+2\,,
\end{align*}
and
\begin{align*}
\tbP_j^{\bfxi}&:=
\tbQ_j^{\bfxi}\,,&
1\leq j\leq2\,,\\
\tbP_j^{\bfxi}&:=
\tbQ_j^{\bfxi}
+\sum_{\ell=1}^{2} \langle \tbQ_j^{\b0};\dD_{\bfxi} \bQ^{\b0}_\ell(\bfxi)\rangle_{L^2}\,\tbQ_\ell^{\bfxi}\,,&
3\leq j\leq r+2\,.
\end{align*}
\end{proof}

From now on we shall enforce the normalization from Lemma~\ref{lem:blocksplit}.

\section{Linear estimates}

We begin the stability analysis with linear estimates, that is, estimates on $(S(t))_{t\geq0}$ the semigroup generated by $L$, in particular proving linear space-modulated stability in the sense made explicit in \cite{R-linKdV}. 

\subsection{Linear phase separation}

With notation from the previous section, we may decompose $S(t)$ according to
\[
S(t)[\bfg]\,=\,(s(t)[\bfg] \cdot \nabla) \ubU\,+\,S_{1}(t)[\bfg] + S_{2}(t)[\bfg]
\]
with
\[
(s(t)[\bfg])(\bfx):= \int_{[-\pi,\pi]^{2}} \frac{\chi(\bfxi)}{\|\bfxi\|}\eD^{\iD \bfx \cdot \bfxi} 
\bp \I_2&\b0_{2,r}\ep
\eD^{t \bDelta_{\bfxi}} 
\bp \left\langle\tbq_{1}^{\bfxi};\widecheck{\bfg}(\bfxi,\cdot) \right\rangle_{L^{2}}\\ 
\vdots\\
\left\langle\tbq_{r+2}^{\bfxi};\widecheck{\bfg}(\bfxi,\cdot) \right\rangle_{L^{2}}  \ep \dD\bfxi\,,
\]
and
\begin{align*}
&(S_{1}(t)[\bfg])(\bfx)\\
&:= \int_{[-\pi,\pi]^{2}} (1-\chi(\bfxi)) \eD^{\iD \bfx \cdot \bfxi} \eD^{t\,L_{\bfxi}} (\widecheck{\bfg}(\bfxi,\cdot))(\bfx) \dD\bfxi 
+\int_{[-\pi,\pi]^{2}}  \chi(\bfxi) \eD^{\iD \bfx \cdot \bfxi} \eD^{t\,L_{\bfxi}} (\I - \Pi_{\bfxi})(\widecheck{\bfg}(\bfxi,\cdot))(\bfx) \dD\bfxi\\[0.5em]
&(S_{2}(t)[\bfg])(\bfx)\\[-1.5em]
&:= \int_{[-\pi,\pi]^{2}} \chi(\bfxi) \eD^{\iD \bfx \cdot \bfxi}  
\bp \frac{\bQ_{1}^{\bfxi}(\bfx) - \bQ_{1}^{\b0}(\bfx)}{\|\bfxi\|}
&\frac{\bQ_{2}^{\bfxi}(\bfx) - \bQ_{2}^{\b0}(\bfx)}{\|\bfxi\|}&\hspace{-0.5em}
\bQ_{3}^{\bfxi}(\bfx)&\hspace{-0.5em}\cdots&\hspace{-0.5em}\bQ_{r+2}^{\bfxi}(\bfx)\ep 
\eD^{t\,\bDelta_{\bfxi}} 
\bp \left\langle \tbq_{1}^{\bfxi};\widecheck{\bfg}(\bfxi,\cdot) \right\rangle_{L^{2}}\\ 
\vdots\\
\left\langle \tbq_{r+2}^{\bfxi};\widecheck{\bfg}(\bfxi,\cdot)\right\rangle_{L^{2}} \ep \dD\bfxi
\end{align*}
where $\chi$ is a smooth function valued in $[0,1]$, compactly supported in a sufficiently small neighborhood of $\b0$ and equal to $1$ in a (smaller) neighborhood of $\b0$.

\subsection{Localized perturbations}

\bpr\label{prop:linestimates}
Assume \cond1-\cond2-\cond3.
\begin{enumerate}
\item There exists $\theta'>0$, such that, for any $(s,s')\in(\R_+)^2$ such that $s'\leq s$, there exists $C_{s',s}$ such that for any $t> 0$
\begin{align*}
\|S_{1}(t)[\bfg]\|_{H^{s}}
&\leq \frac{C_{s',s}}{(\min(\{1,t\}))^{\frac{(s-s')}{2}}}\,\eD^{-\theta'\,t}
\|\bfg\|_{H^{s'}}\,.
\end{align*}
\item For any $s\in\R_+$, any $\beta\in\N^2$ and any $\gamma\in\N^2$, there exists $C_{s,\beta,\gamma}$ such that for any $1\leq q \leq 2 \leq p \leq +\infty$, and any $t\geq0$
\begin{align*}
\|S_{2}(t)[\d_\bfx^\beta\,\cJ^\gamma\,\bfg]\|_{W^{s,p}}
&\leq \frac{C_{s,\beta,\gamma}}{(1+t)^{\frac{|\gamma|}{2}+\frac{1}{q}-\frac{1}{p}}}\,\|\bfg\|_{L^{q}}\,,
\end{align*}
and if $|\beta|>0$ or $\bfg$ vanishes identically in its first $r$ components
\begin{align*}
\|S_{2}(t)[\d_\bfx^\beta\,\cJ^\gamma\bfg]\|_{W^{s,p}}
&\leq \frac{C_{s,\beta,\gamma}}{(1+t)^{\frac{1+|\gamma|}{2}+\frac{1}{q}-\frac{1}{p}}}\,\|\bfg\|_{L^{q}}\,.
\end{align*}
\item For any $\alpha\in\N^2$, any $\beta\in\N^2$, any $\gamma\in\N^2$ and any $\ell\in\N$, there exists $C_{\alpha,\ell,\beta,\gamma}$ such that for any $1\leq q \leq 2 \leq p \leq +\infty$, and any $t\geq0$
\begin{enumerate}
\item provided that $|\alpha|+|\gamma|+\ell\geq 1$ or $1/q-1/p>1/2$ 
\begin{align*}
\|\,\d_\bfx^\alpha\,\d_t^\ell\,s(t)[\d_\bfx^\beta\,\cJ^\gamma\,\bfg]\|_{L^p}
&\leq \frac{C_{\alpha,\ell,\beta,\gamma}}{(1+t)^{\frac{|\alpha|+|\gamma|+\ell-1}{2}+\frac{1}{q}-\frac{1}{p}}}\,\|\bfg\|_{L^{q}}\,,
\end{align*}
\item if $|\beta|>0$ or $\bfg$ vanishes identically in its first $r$ components
\begin{align*}
\|\,\d_\bfx^\alpha\,\d_t^\ell\,s(t)[\d_\bfx^\beta\,\cJ^\gamma\bfg]\|_{L^p}
&\leq \frac{C_{\alpha,\ell,\beta,\gamma}}{(1+t)^{\frac{|\alpha|+|\gamma|+\ell}{2}+\frac{1}{q}-\frac{1}{p}}}\,\|\bfg\|_{L^{q}}\,.
\end{align*}
\end{enumerate}
\end{enumerate}
\epr

We recall that $\cJ$ commutes with $\d_\bfx$, $\d_t$, $s(t)$, $S_1(t)$, $S_2(t)$, \emph{etc.}

Compared to \cite[Proposition~2.2]{MR} algebraic decay estimates for a general $\bfg$ are deteriorated by a $(1+t)^{-1/2}$ factor. Compared to \cite[Proposition~3.3]{JNRZ-conservation}, there are two main differences:
\begin{itemize}
\item the slope of $1/q-1/p$ depends on the dimension, being equal to $d/2$ in dimension $d$
\item the restriction $1\leq q \leq 2 \leq p \leq +\infty$ is not a pure feature of the chosen proof, but a hard restriction that cannot be relaxed at need without deteriorating the heat-type decay rates. 
\end{itemize}

\begin{proof}
To prove the first point, it is sufficient to combine an $H^{s'}\to H^s$ bound for $0<t\leq 1$ with an $H^s\to H^s$ bound for $t\geq 0$. Moreover the former follows from the parabolicity of $L$ (combined with bounds on $s(t)$ and $S_2(t)$ proved below). In turn, the latter may be derived, through Parseval's identity and exponential bounds from Proposition~\ref{p:Bloch-bounds}.

All the algebraic decay rates arise from the fact when $1\leq r\leq\infty$ and ($\eta>-2/r$ or $\eta=0$)
\[
\left\|\bfxi\mapsto \|\bfxi\|^{\eta}\,\eD^{-\theta t \|\bfxi\|^{2}}\right\|_{L^{r}_\bfxi}
\lesssim (1+t)^{-\left(\frac{\eta}{2}+\frac{1}{r}\right)}\,.
\]
For instance, to prove the first part of the second point, by integration by parts in scalar products, from \cond3, Hausdorff-Young and H\"older inequalities one derives 
\begin{align*}
\|S_{2}(t) [\d_\bfx^\beta\,\cJ^\gamma\bfg]  \|_{W^{s,p}} 
&\lesssim \left\|\bfxi\mapsto \eD^{-\theta t \|\bfxi\|^{2}} \|\bfxi\|^{|\gamma|}\|\widecheck{\bfg}(\bfxi,\cdot)\|_{L^{q}}\right\|_{L^{p'}_\bfxi}\\
&\lesssim  \left\|\bfxi\mapsto \|\bfxi\|^{|\gamma|}\,\eD^{-\theta t \|\bfxi\|^{2}}\right\|_{L^{r}_\bfxi}\ \times\ \|\widecheck{\bfg} \|_{L^{q'}_{\bfxi} L^{q}_{\bfx}}
\lesssim (1+t)^{-\left(\frac{|\gamma|}{2}+\frac{1}{r}\right)}\,\|\bfg\|_{L^q}
\end{align*}
with $p'$, $q'$ Lebesgue conjugate respectively to $p$ and $q$, and $1/r=1/p'-1/q'=1/q-1/p$. Hence the second bound. The second part of the second point follows from a similar bound benefiting from an extra $\|\bfxi\|$ factor due to the cancellation arising from the extra $\bfg$ structure and $(\tbq_{1}^{\b0},\cdots,\tbq_{r+2}^{\b0})\equiv (\b0,\b0,\beD_1,\cdots,\beD_{r})$. The third point is proved similarly starting from
\begin{align*}
(\d_\bfx^\alpha\,\d_t^\ell s(t)&[\d_\bfx^\beta\,\cJ^\gamma\,\bfg])(\bfx)\\
&
=\int_{[-\pi,\pi]^{2}} \frac{\chi(\bfxi)}{\|\bfxi\|}\eD^{\iD \bfx \cdot \bfxi}(\iD\bfxi)^\alpha\,\bfxi^\beta\,
\bp \I_2&\b0_{2,r}\ep
(\bDelta_\bfxi)^\ell\,
\eD^{t \bDelta_{\bfxi}} 
\bp \left\langle\tbq_{1}^{\bfxi};(\d_\bfx+\iD\bfxi)^\beta\widecheck{\bfg}(\bfxi,\cdot)\right\rangle_{L^{2}}\\ 
\vdots\\
\left\langle\tbq_{r+2}^{\bfxi};(\d_\bfx+\iD\bfxi)^\beta\widecheck{\bfg}(\bfxi,\cdot)\right\rangle_{L^{2}}  \ep \dD\bfxi\,.
\end{align*}
\end{proof}

\subsection{Perturbation by phase modulation}

We restrict here to estimates useful to analyze the critical-decay case.

Throughout we implicitly assume that $\bfphi$ has no affine component at $\infty$, in the sense that $\bfphi=\Delta^{-1}\Delta\bfphi$. Consistently, the phases built with $s(t)$ also satisfy the latter condition.

\bpr\label{prop:linphasedecay}
Assume \cond1-\cond2-\cond3.
\begin{enumerate}
\item There exists $\theta'>0$, such that, for any $(s,s')\in\R_+$ such that $s'+1\leq s$ and any $(p_0,p_1)$ such that $2<p_0<p_1<\infty$, there exists $C_{p_0,p_1,s,s'}$ such that for any $t>0$, and any $p_0\leq p\leq p_1$,
\begin{align*}
\|S_{1}(t)[(\bfphi \cdot \nabla) \ubU]\|_{W^{s,p}}
&\leq \frac{C_{p_0,p_1,s,s'}}{(\min(\{1,t\}))^{\frac{(s-(s'+1))}{2}}}\,\eD^{-\theta'\,t}
\,\|\Delta\bfphi\|_{L^1\cap H^{s'}}\,.
\end{align*}
\item For any $\alpha\in\N^2$, any $s\in\R_+$ and any $2 \leq p \leq +\infty$ such that $|\alpha|+1-\tfrac2p>0$, there exists $C_{p,\alpha,s}$ such that for any $t\geq0$
\begin{align*}
\|\cJ^\alpha S_{2}(t)[(\bfphi \cdot \nabla) \ubU]\|_{W^{s,p}}
&\leq \frac{C_{p,\alpha,s}}{(1+t)^{\frac{|\alpha|+1}{2}-\frac{1}{p}}}\,\|\Delta\bfphi\|_{L^1}\,.
\end{align*}
\item For any $\alpha\in\N^2$, any $\beta\in\N^2$, any $\ell\in\N$ and any $2 \leq p \leq +\infty$ such that $|\alpha|+|\beta|+\ell-\tfrac2p>0$, there exists $C_{p,\alpha,\beta,\ell}$ such that for any $t\geq0$
\begin{align*}
\|\,\d_\bfx^\alpha\,\d_t^\ell\,\cJ^\beta\,s(t)[(\bfphi \cdot \nabla) \ubU]\|_{L^p}
&\leq \frac{C_{p,\alpha,\ell}}{(1+t)^{\frac{|\alpha|+|\beta|+\ell}{2}-\frac{1}{p}}}\,\|\Delta\bfphi\|_{L^1}\,.
\end{align*}
\end{enumerate}
\epr

To ease comparisons with bounds of Proposition~\ref{prop:linestimates}, we point out that $\|\Delta\bfphi\|_{L^1\cap H^{s}}$ should be thought as a relaxed version of $\|\nabla\bfphi\|_{H^{s+1}}$. Note moreover that the condition $|\alpha|+|\beta|+\ell-\tfrac2p>0$ may be written more explicitly as $|\alpha|+|\beta|+\ell\geq 2$ or ($|\alpha|+|\beta|+\ell=1$ and $p>2$).

\begin{proof}
To establish various bounds it is convenient to single out the low-frequency part of $\bfphi$, according to
\begin{align*}
\bfphi&\,=\,\bfphi_{LF}+\bfphi_{HF}\,,&
\widehat{(\bfphi_{LF})}&=\chi \widehat{\bfphi}\,.
\end{align*} 
The contribution of $\bfphi_{HF}$ to the first bound may be deduced from the corresponding estimate in Proposition~\ref{prop:linestimates}. Indeed, since $2\leq p<\infty$,
\begin{align*}
\|S_{1}(t)[(\bfphi_{HF} \cdot \nabla) \ubU]\|_{W^{s,p}}
&\lesssim \|S_{1}(t)[(\bfphi_{HF} \cdot \nabla) \ubU]\|_{H^{s+1}}\,,\\
\|(\bfphi_{HF} \cdot \nabla) \ubU\|_{H^{s'+2}}
&\lesssim\|\bfphi_{HF}\|_{H^{s'+2}}
\lesssim\|\Delta\bfphi\|_{H^{s'}}\,.
\end{align*}
The analysis of the contribution of $\bfphi_{LF}$ requires more care. To begin with, we recall that 
\begin{align}\label{eq:low-phi}
\widecheck{((\bfphi_{LF} \cdot \nabla) \ubU)}(\bfxi,\bfx) = 
(\widehat{(\bfphi_{LF})}(\bfxi) \cdot \nabla) \ubU(\bfx)
\end{align}
and observe that this may be used to gain an extra $\|\bfxi\|$-factor in the second part of the definition of $S_1$ through
\[
(\I - \Pi_{\bfxi})(\widecheck{((\bfphi_{LF} \cdot \nabla) \ubU)}(\bfxi,\cdot))
\,=\,(\I - \Pi_{\bfxi})(\Pi_{\b0} - \Pi_{\bfxi})
(\widecheck{((\bfphi_{LF} \cdot \nabla) \ubU)}(\bfxi,\cdot))
\] 
(which holds since $\d_1\ubU$ and $\d_2\ubU$ lie in the range of $\Pi_0$ and $(\I - \Pi_{\bfxi})$ is a projector). Moreover, an extra $\|\bfxi\|$-factor in the first part of $S_1$ is readily obtained from the trivial $(1-\chi(\bfxi))\lesssim \|\bfxi\|$. With this in hands, from Hausdorff-Young inequalities and the embedding $H^{s+1}_{\rm per}\hookrightarrow W^{s,p}_{\rm per}$, one derives
\begin{align*}
\|S_{1}(t)[(\bfphi_{LF} \cdot \nabla) \ubU]\|_{W^{s,p}}
\lesssim
\eD^{-\theta'\,t}\,\,\|\Delta\bfphi_{LF}\|_{L^1}\,
\left\|\bfxi\mapsto \|\bfxi\|^{-1}\right\|_{L^{p'}_\bfxi}
\lesssim
\eD^{-\theta'\,t}\,\,\|\Delta\bfphi\|_{L^1}
\end{align*}
since $p>2$, thus $p'<2$. This achieves the proof of the first bound.

The contribution of $\bfphi_{HF}$ to the second bound may also be deduced from the corresponding estimate in Proposition~\ref{prop:linestimates}. Indeed
\begin{align}\label{estim:high-phi}
\|(\bfphi_{HF} \cdot \nabla) \ubU\|_{L^1}
&\lesssim\|\bfphi_{HF}\|_{L^1}
\lesssim\|\Delta\bfphi\|_{L^1}\,.
\end{align}
By using the cancellations brought by \eqref{eq:low-phi} the analysis in terms of $\|\Delta\bfphi\|_{L^1}$ of the contribution of $\bfphi_{LF}$ to the second bound becomes similar to the one in terms of $\|\cJ^\alpha\bfg\|_{L^1}$ of the contribution of a general $\bfg$ to $s(t)[\cJ^\alpha\bfg]$ in Proposition~\ref{prop:linestimates}.

The third bound is proved similarly with an extra $\|\bfxi\|^{|\beta|+\ell-1}$ factor. 
\end{proof}

\subsection{Short-time phase estimate}

The final linear estimate we need is a short-time bound on $s(t)$. Let us anticipate on the nonlinear scheme and mention that we only need it insofar as to enforce $\phi(0,\cdot)=\phi_0$ in the Duhamel formulation.

\bl\label{l:boundarylayer}
Assume \cond1-\cond2-\cond3. For any $\alpha\in\N^2$, any $\ell\in\N$ and any $2 \leq p \leq +\infty$ such that $p>2$ if $|\alpha|=0$ and $p<\infty$ if $|\alpha|\geq2$, there exists $C_{p,\alpha}$ such that for any $t\geq 0$
\begin{align*}
\|\,\d_\bfx^\alpha\,\left(s(t)[(\bfphi \cdot \nabla)\ubU]-\bfphi\right)\|_{L^p}
&\leq C_{p,\alpha}\,\|\Delta\bfphi\|_{L^1\cap W^{(|\alpha|-2)_+,p}}
\,\begin{cases}
\,(1+t)^{\frac12\,\left(1-\left(|\alpha|-\frac2p\right)\right)_+}&\quad\textrm{if } |\alpha|-\tfrac2p\neq 1\\
\ln(2+t)&\quad\textrm{otherwise }
\end{cases}\,.
\end{align*}
\el

\begin{proof}
We first take a detour and show that it suffices to bound $s(0)[(\bfphi\cdot\nabla)\ubU ]-\bfphi$.
To do so, we integrate the bounds on $\d_\bfx^\alpha\d_t\,s$ from Proposition~\ref{prop:linphasedecay} and are led to
\begin{align*}
\|\,\d_\bfx^\alpha\,(s(t)-s(0))[(\bfphi \cdot \nabla)\ubU]\|_{L^p}
&\lesssim 
\|\Delta\bfphi\|_{L^1}\,\int_0^t\frac{\dD\tau}{(1+\tau)^{\frac{|\alpha|+1}{2}-\frac{1}{p}}}\,\,,
\end{align*}
which predicts the growth time rates.

To bound $s(0)[(\bfphi\cdot\nabla)\ubU ]-\bfphi$, we split the phase $\phi$ into low and high frequency components as in the proof of Proposition~\ref{prop:linphasedecay}. For the high-frequency contribution, we rely on the triangle inequality and observe that the proof of Proposition~\ref{prop:linphasedecay} yields
\begin{align*}
\|\,\d_\bfx^\alpha\,s(t)[(\bfphi_{HF} \cdot \nabla)\ubU]\|_{L^p}
&\lesssim 
\frac{\|\Delta\bfphi\|_{L^1}}{(1+t)^{\frac{|\alpha|}{2}+1-\frac{1}{p}}}
\lesssim 
\,\|\Delta\bfphi\|_{L^1}\,\,,
\end{align*}
whilst the conditions on $p$ ensure
\[
\|\bfphi_{HF}\|_{W^{|\alpha|,p}}
\lesssim \|\Delta\bfphi\|_{L^1\cap W^{(|\alpha|-2)_+,p}}\,,
\]

For the low-frequency contribution of $\bfphi$, we observe that
\begin{align*}
\left(s(0)[(\bfphi_{LF} \cdot \nabla) \ubU] - \bfphi_{LF}\right)(\bfx)
&=\int_{[-\pi,\pi]^{2}} \eD^{\iD \bfxi \cdot\bfx}\,(\chi(\bfxi)-1)\,\widehat{(\bfphi_{LF})}(\bfxi)\,\dD\bfxi\\
&\quad+\int_{[-\pi,\pi]^{2}}\eD^{\iD \bfxi \cdot\bfx}\,\chi(\bfxi) \,
\bp \left\langle \tbq_{1}^{\bfxi} - \tbq_{1}^{\b0};(\widehat{(\bfphi_{LF})}(\bfxi) \cdot \nabla) \ubU\right\rangle_{L^{2}_{per}} \\ 
\left\langle \tbq_{2}^{\bfxi} - \tbq_{2}^{\b0} ;(\widehat{(\bfphi_{LF})}(\bfxi) \cdot \nabla) \ubU\right\rangle_{L^{2}_{per}} \ep 
\dD\,\bfxi\,.
\end{align*}
Thus
\[
\|\,\d_\bfx^\alpha\,\left(s(0)[(\bfphi \cdot \nabla)\ubU]-\bfphi\right)\|_{L^p}
\lesssim \|\Delta\bfphi\|_{L^1}
\,\times\,\left\|\bfxi\mapsto \|\bfxi\|^{|\alpha|-1}\right\|_{L^{p'}_\bfxi}\,.
\]
Hence the result (since $p'<2$ when $|\alpha|=0$).
\end{proof}

\section{Nonlinear stability}

\subsection{Nonlinear separation of the phase}

Our nonlinear analysis begins by reformulating Equation~\eqref{eq:rd} in terms of $\bfphi$ and $\bV$ such that $\bW(t,\cdot-\bfphi(t,\cdot))=\ubU+\bV(t,\cdot)$.

At the beginning the algebraic manipulations are identical to the ones of \cite[Section~2.2]{MR}. Namely, we introduce
\begin{align*}
\cA[\bW]
&=\cA(\bW,\nabla\bW,\nabla^2\bW)
\,:=\,
\transp{\nabla}\left(\transp{\ubK}\ubK\,\nabla\bW\right)
+\transp{\nabla}\left(\transp{\ubK}\,
(\bG(\bW)+\ubc\,\transp{\bW})\right)+\bff(\bW)
\end{align*}
and consider its image under a change of variable $\bPhi$
\begin{align}\label{def:bA}
\bA[\tbW,\bPhi]
&=\bA(\tbW,\nabla\tbW,\nabla^2\tbW,\nabla\bPhi,\nabla^2\bPhi)
:=(\cA[\tbW\circ \bPhi^{-1}])\circ\bPhi\\\nn
&\,=\,
|\nabla\bPhi|^{-1}
\transp{\nabla}\left(|\nabla\bPhi|\,\transp{(\ubK\,[\nabla\bPhi]^{-1})}(\ubK\,[\nabla\bPhi]^{-1})\nabla\tbW\right)\\\nn
&\quad+\ |\nabla\bPhi|^{-1}
\transp{\nabla}\left(|\nabla\bPhi|\,\transp{(\ubK\,[\nabla\bPhi]^{-1})}
(\bG(\tbW)+\ubc\,\transp{\tbW})\right)+\bff(\tbW)\,.
\end{align}
At the linear level, the key observation is that
\[
-\bL_{\bPhi}\bA[\ubU,\Id](\bfphi)
\,=\,\bL_\bW\cA [\ubU]((\bfphi\cdot\nabla)\ubU)-(\bfphi\cdot\nabla)(\cA[\ubU])
\] 
thus, since $\cA[\ubU]\equiv 0$,
\begin{equation}\label{eq:key-lin}
\bL_{(\tbW,\bPhi)}\bA[\ubU,\Id](\bV,-\bfphi)
\,=\,\bL_\bW \cA [\ubU](\bV+(\bfphi\cdot\nabla)\ubU)
\,=\,L(\bV+(\bfphi\cdot\nabla)\ubU)\,,
\end{equation}
where $\bL_{\bPhi}$, $\bL_\bW$ and $\bL_{(\tbW,\bPhi)}$ stand for linearization operators. With this in hands, we may rephrase \eqref{eq:rd}.

\bl\label{lem:cancellation-separation}
Let $\bW$ and $(\bfphi,\bV)$ be smooth functions such that
\be\label{pertvar} 
\bW(t,\bfx-\bfphi(t,\bfx))\ =\ \ubU(\bfx)+\,\bV(t,\bfx)\,,
\ee
and for any $t$, $\|\nabla\bfphi(t,\cdot)\|_{L^\infty(\R^2)}<1$. Then $\bW$ satisfies \eqref{eq:rd} if and only if $(\bfphi,\bV)$ satisfies
\begin{equation}\label{veq}
\left(\d_t-L\right)(\bV+(\bfphi \cdot \nabla) \ubU)=\cN[\bV,\bfphi]\,,
\end{equation}
or equivalently
\begin{equation}\label{veq0}
\d_t\bV-L_{\bfphi}\bV
\,=\,-(\bfphi_t \cdot [\I_2-\nabla\bfphi]^{-1}\nabla)(\ubU+\bV)
+\bA[\ubU,\Id-\bfphi]-\bA[\ubU,\Id]+\cN_0[\bV,\bfphi]\,,
\end{equation}
with 
\begin{align*}
L_{\bfphi}\bV
&:=\bL_{\tbW}\bA[\ubU,\Id-\bfphi](\bV)\,,\\
\cN_0[\bV,\bfphi]&=\cN_0(\nabla\bfphi,\nabla^2\bfphi,\bV,\nabla\bV)
:=\bA[\ubU+\bV,\Id-\bfphi]
-\bA[\ubU,\Id-\bfphi]-\bL_{\tbW}\bA[\ubU,\Id-\bfphi](\bV)\,,\\
\cN[\bV,\bfphi]&=\cN(\bfphi_t,\nabla\bfphi,\nabla^2\bfphi,\bV,\nabla\bV,\nabla^2\bV)\\
&:=\ \left(\bL_{\tbW}\bA[\ubU,\Id-\bfphi]-\bL_{\tbW}\bA[\ubU,\Id]\right)(\bV)
+\cN_0[\bV,\bfphi]
-(\bfphi_t \cdot \nabla\bfphi\,[\I_2-\nabla\bfphi]^{-1}\nabla)\ubU\\[0.5em]
&
\quad-(\bfphi_t\cdot [\I_2-\nabla\bfphi]^{-1}\nabla) \bV
+\bA[\ubU,\Id-\bfphi]-\bA[\ubU,\Id]-\bL_{\bPhi}\bA[\ubU,\Id](-\bfphi)\,.
\end{align*}
\el

\begin{proof}
To begin with, note that if $\bPsi$ is defined by $\bPsi(t,\cdot):=\bPhi(t,\cdot)^{-1}$, then
\begin{align*}
\nabla\bPsi(t,\bfx)&=[\nabla\bPhi(t,\bPsi(t,\bfx))]^{-1}\,,&
\d_t\bPsi(t,\bfx)&=-[\transp{(\nabla\bPhi(t,\bPsi(t,\bfx)))}]^{-1}
\,\d_t\bPhi(t,\bPsi(t,\bfx))\,.
\end{align*}
Define $\tbW$ by $\bW(t,\bPhi(t,\bfx))=\tbW(t,\bfx)$ or equivalently $\tbW(t,\bfx)=\bW(t,\bPsi(t,\bfx))$. Then $\bW$ solves \eqref{eq:rd} if and only if $(\bPhi,\tbW)$ satisfies
\begin{align*}
\tbW_t- (\bPhi_t \cdot [\nabla\bPhi]^{-1}\nabla) \tbW
&\ =\ \bA[\tbW,\bPhi]\,.
\end{align*}
Inserting $\tbW=\ubU+\bV$ and $\bPhi=\Id-\bfphi$, we readily deduce \eqref{veq0} and derive \eqref{veq} by combining it with \eqref{eq:key-lin}.
\end{proof}

The form \eqref{veq} is adapted to the large-time analysis whereas the form \eqref{veq0} is used in nonlinear regularity estimates. A key consequence of Lemma~\ref{lem:cancellation-separation} is that, as long as $\|\nabla\bfphi\|_{L^\infty}<1$, \eqref{eq:rd} is equivalently written as $\bV(0,\cdot)=\bV_0$, $\bfphi(0,\cdot)=\bfphi_0$ and
\[
\bV(t,\cdot)+(\bfphi(t,\cdot) \cdot \nabla) \ubU
\,=\,S(t)[\bV_0+(\bfphi_0 \cdot \nabla) \ubU]
+\int_0^tS(t-\tau)\,\cN[\bV(\tau,\cdot),\bfphi(\tau,\cdot)]\,\dD \tau\,.
\]
At this stage,  we would like to simply use the semigroup splitting of the linear analysis so as to split the foregoing nonlinear equation while enforcing $\bV(0,\cdot)=\bV_0$, $\bfphi(0,\cdot)=\bfphi_0$. To do so, we pick $\tchi$ a smooth function on $\R_+$ valued in $[0,1]$, compactly supported in $[0,1]$ and equal to $1$ on $[0,\tfrac12]$. Then, we consider
\begin{align}\label{bfphi}
\bfphi(t,\cdot)&=s(t)[\bV_{0} + (\bfphi_{0} \cdot \nabla) \ubU] 
+ \int_{0}^{t} s(t-\tau)\cN[\bV(\tau,\cdot),\bfphi(\tau,\cdot)]\dD \tau\\
&\ +\tchi(t)\,\left(\bfphi_{0}-s(t)[\bV_{0} + (\bfphi_{0} \cdot \nabla) \ubU]\right)
\nn\\\label{V}
\bV(t,\cdot)&=(S_{1}+S_{2})(t) \left[ \bV_{0}+(\bfphi_{0} \cdot \nabla) \ubU \right] + \int_{0}^{t} (S_{1}+S_{2})(t-\tau)\cN[\bV(\tau,\cdot),\bfphi(\tau,\cdot)]\dD\tau\\
&\ -\tchi(t)\,\Big(\left(\bfphi_{0}-s(t)[\bV_{0} + (\bfphi_{0} \cdot \nabla) \ubU]\right)\cdot\nabla\Big)\ubU
\nn
\end{align}
and observe that, as long as $\|\nabla\bfphi\|_{L^\infty}<1$, \eqref{bfphi}-\eqref{V} imply that $\bW$ defined by 
\[
\bW(t,\cdot):=(\ubU+\bV(t,\cdot))\circ (\Id-\bfphi(t,\cdot))^{-1}
\]
satisfies \eqref{eq:rd} with $\bW(0,\cdot):=(\ubU+\bV_0)\circ (\Id-\bfphi_0)^{-1}$. 

Now we would like to use \eqref{veq0} to derive a nonlinear high-frequency damping estimate and insert the linear estimates of the previous section in \eqref{bfphi}-\eqref{V}. Yet, to make the most of the latter, we must add to the computations of \cite{MR} a further study of the conservative structure of $\cN[\bV,\bfphi]$, reminiscent of arguments from \cite{JNRZ-conservation}. To do so we begin by splitting $\bA$ into flux, source and commutator terms 
\[
\bA[\tbW,\bPhi]\,=\,\transp{\nabla}\left(\bA_f[\tbW,\bPhi]\right)+\bA_s[\tbW]+\bA_c[\tbW,\bPhi]
\]
where
\begin{align*}
\bA_f[\tbW,\bPhi]
&=\bA_f(\tbW,\nabla\tbW,\nabla\bPhi)
\,:=\,
\transp{(\ubK\,[\nabla\bPhi]^{-1})}(\ubK\,[\nabla\bPhi]^{-1})\nabla\tbW
+\transp{(\ubK\,[\nabla\bPhi]^{-1})}
(\bG(\tbW)+\ubc\,\transp{\tbW})\\
\bA_c[\tbW,\bPhi]
&=\bA(\tbW,\nabla\tbW,\nabla\bPhi,\nabla^2\bPhi)\\
&\,:=\,
\left[|\nabla\bPhi|^{-1},
\transp{\nabla}\right]\left(|\nabla\bPhi|\,\transp{(\ubK\,[\nabla\bPhi]^{-1})}
\left((\ubK\,[\nabla\bPhi]^{-1})\nabla\tbW+\bG(\tbW)+\ubc\,\transp{\tbW}\right)\right)\,,
\end{align*}
and $\bA_s[\tbW]:=\bff(\tbW)$. 

The following lemma provides the required structure. Note that by design flux-type terms are conservative, source-type terms vanish identically in their first $r$ components and commutator terms contain second-order derivatives of phase shifts as linear factors.

\bl\label{lem:structure}
With assumptions and notation from Lemma~\ref{lem:cancellation-separation},
\[
\cN[\bV,\bfphi]\,=\,\transp{\nabla}\left(\cN_f[\bV,\bfphi]\right)+\cN_s[\bV]+\cN_c[\bV,\bfphi]
\]
with
\begin{align*}
\cN_f[\bV,\bfphi]&=\cN_f(\bfphi_t,\nabla\bfphi,\bV,\nabla\bV)\\
&:=\ \left(\bL_{\tbW}\bA_f[\ubU,\Id-\bfphi]-\bL_{\tbW}\bA_f[\ubU,\Id]\right)(\bV)
-\ubU\,\transp{(\bfphi_t)}\,\nabla\bfphi\,[\I_2-\nabla\bfphi]^{-1}\\[0.5em]
&
\quad-\bV\,\transp{(\bfphi_t)}[\I_2-\nabla\bfphi]^{-1} 
+\bA_f[\ubU,\Id-\bfphi]-\bA_f[\ubU,\Id]-\bL_{\bPhi}\bA_f[\ubU,\Id](-\bfphi)\\[0.5em]
&\quad+\bA_f[\ubU+\bV,\Id-\bfphi]
-\bA_f[\ubU,\Id-\bfphi]-\bL_{\tbW}\bA_f[\ubU,\Id-\bfphi](\bV)\,,\\[0.5em]
\cN_s[\bV]&=\cN_s(\bV)
:=\bA_s[\ubU+\bV]-\bA_s[\ubU]-\bL_{\tbW}\bA_s[\ubU](\bV)\,,\\[0.5em]
\cN_c[\bV,\bfphi]&=\cN(\bfphi_t,\nabla\bfphi_t,\nabla\bfphi,\nabla^2\bfphi,\bV,\nabla\bV)\\
&:=\ \left(\bL_{\tbW}\bA_c[\ubU,\Id-\bfphi]-\bL_{\tbW}\bA_c[\ubU,\Id]\right)(\bV)
+\ubU\,\transp{\nabla}\left(\transp{(\bfphi_t)}\,\nabla\bfphi\,[\I_2-\nabla\bfphi]^{-1}\right)\\[0.5em]
&
\quad+\bV\,\transp{\nabla}\left(\transp{(\bfphi_t)}[\I_2-\nabla\bfphi]^{-1}\right)
+\bA_c[\ubU,\Id-\bfphi]-\bA_c[\ubU,\Id]-\bL_{\bPhi}\bA_c[\ubU,\Id](-\bfphi)\\[0.5em]
&\quad+\bA_c[\ubU+\bV,\Id-\bfphi]
-\bA_c[\ubU,\Id-\bfphi]-\bL_{\tbW}\bA_c[\ubU,\Id-\bfphi](\bV)\,.
\end{align*}
\el

For later use note that when $\|\nabla\bfphi\|$ is small and bounded away from $1$ and $\bV$ varies in a compact, one has the pointwise bounds
\begin{align*}
\| \cN_f[\bV,\bfphi] \| &\lesssim \| \bfphi_{t} \| \left(\|\bV\|+\| \nabla \bfphi \|\right)+\|\nabla \bV \|\| \nabla \bfphi\|+\left(\|\bV\|+\| \nabla \bfphi \|\right)^2\,,\\
\| \cN_s[\bV,\bfphi] \| &\lesssim \|\bV\|^2\,,\\
\| \cN_c[\bV,\bfphi] \| &\lesssim \|\nabla\bfphi_{t} \| \left(\|\bV\|+\| \nabla \bfphi \|\right)+
\|\nabla^2\bfphi\| \left(\|\bV\|+\|\nabla \bV \|+\| \nabla \bfphi \|+\| \bfphi_{t} \| \right)\,.
\end{align*}

\subsection{Dissipation of nonlinear high-frequency damping estimates}

To close estimates in regularity we shall perform energy estimates on \eqref{veq0}. Note however that as in \cite{MR} but unlike \cite{JNRZ-conservation}, the localization we may use to do so is limited. Indeed as apparent in Proposition~\ref{prop:linphasedecay}, from $\Delta \bfphi\in L^1$ one may only obtain $\nabla\bfphi\in L^p$ when $p>2$. This limitation is transferred to $\bV$ and its derivatives.

Thus we must perform the energy estimates in $L^p$, $p>2$. Since energy estimates are cleaner for integer powers of $2$, the choice we have made in Theorem~\ref{theo:critmain} is to rely on $L^4$ estimates, as in \cite{MR}. 

Since $\d_t(\|\bV\|)^4=4\,\|\bV\|^2\bV\cdot\d_t\bV$, the key ingredient in deriving nonlinear high-frequency damping estimates is provided by the following lemma.

\bl\label{l:energy-estimate}
\begin{enumerate}
\item For any $0\leq\eta_0<1$ and any $\ell\in\N$, there exist $\theta>0$ and $C\geq 0$ such that for any $\bfphi$ such that $\|\nabla\bfphi\|_{L^\infty(\R^2;\cM_2(\R))}\leq \eta_0$ and any $\bV\in W^{\ell,4}(\R^2;\R^n)$ such that $\cD_\ell(\bV)<+\infty$ and $L_{\bfphi}\bV\in W^{\ell,4}(\R^2;\R^n)$
\bas
\sum_{|\alpha|=\ell} \int_{\R^2}\,\|\d^\alpha\bV\|^2\,\d^\alpha\bV\cdot \d^\alpha(L_{\bfphi}\bV)
&\leq -\theta\,\cD_\ell(\bV)^4 \,+\,C\, \|\bV\|_{L^4}^4\,\left(1+\|\nabla^2\bfphi\|_{L^\infty}^{2\,(2\ell+1)}\right)\\
&+C \|\nabla\bfphi\|_{W^{\ell,4}}^4\,\|\bV\|_{W^{1,\infty}}^{4}\,
\eas
where
\[
\cD_\ell(\bV):=\left(\sum_{|\alpha|=\ell}\,\int_{\R^2}\,\|\d^\alpha\bV\|^2\,\|\nabla\d^\alpha\bV\|^2\,\right)^{\frac14}\,.
\]
\item For any $\ell\in\N$, there exists $C$ such that for any $\bV\in W^{\ell+1,4}(\R^2;\R^n)$ and any $j\in\N$, $j\leq \ell$,
\begin{align*}
\|\nabla^j\bV\|_{L^4}
&\leq C\,\|\bV\|_{L^4}^{1-\alpha}
\,\cD_\ell(\bV)^{\alpha}\,,&
\textrm{with }\alpha:=\frac{j}{\ell+\tfrac12}\,.
\end{align*}
\end{enumerate}
\el

Lemma~\ref{l:energy-estimate} is identical to \cite[Lemma~2.8]{MR} and the proof provided there also applies to our present case. We omit to repeat it here.

The content of the lemma is that up to lower order terms $W^{\ell,4}$-estimates provide a dissipation similar, from the point of view of Sobolev embeddings, to a $W^{\ell+\frac12,4}$ damping. Combined with \eqref{veq0} and the low-frequency nature it implies that $\|\bV(t,\cdot)\|_{W^{2,4}}$ decays at least as fast as $\|(\nabla\bfphi,\bfphi_t,\bV)(t,\cdot)\|_{L^{4}}$ does.

\subsection{Proof of Theorem~\ref{theo:critmain}}\label{s:critmain}

To keep the discussion focused on the large-time analysis we shall take the local existence of solutions to \eqref{bfphi}-\eqref{V} for granted and prove Theorem~\ref{theo:critmain} by a continuity argument. We stress however that one may readily extract from \cite[Subsection~3.2]{MR} the needed local existence result with blow-up criterion expressed in terms of $\|(\bV,\nabla \bfphi,\bfphi_t)\|_{W^{1,\infty}}$. Since $W^{2,4}$ embeds continuously in $W^{1,\infty}$, proving a $W^{2,4}$ control is sufficient to guarantee global existence.

Before digging into the proof, we briefly discuss its spirit. Our strategy is to close $W^{2,4}$ bounds alone by a continuity argument then to use those to prove the remaining $L^p$ bounds. When doing so, the propagation of algebraic decay rates relies on the basic 
\begin{align*}
\int_0^{t/2}\frac{1}{(1+t-\tau)^\alpha}\frac{1}{(1+\tau)^\beta}\dD\tau &\lesssim\frac{1}{(1+t)^{\alpha+\beta-1}}\,,& 0\leq \alpha\,,\quad 0\leq \beta<1\,,\\
\int_{t/2}^t\frac{1}{(1+t-\tau)^\alpha}\frac{1}{(1+\tau)^\beta}\dD\tau &\lesssim\frac{1}{(1+t)^{\alpha+\beta-1}}\,,& 0\leq \alpha<1\,,\quad 0\leq \beta\,.
\end{align*}
When estimating critical terms any other estimate would fail to close in decay and thus one must be able to distribute decay so as to enforce the respective constraints $\beta<1$ and $\alpha<1$. For flux-type and source-type nonlinear terms the final count is similar to the nonlinear analysis of \cite[Subsection~3.2]{MR} and follows from an $L^2\to L^4$ bound resulting in the special case 
\[
\int_0^t\frac{1}{(1+t-\tau)^{\frac34}}\frac{1}{(1+\tau)^{\frac12}}\dD\tau 
\lesssim\frac{1}{(1+t)^{\frac14}}\,.
\]
However the case of commutator-type terms is significantly different. Since they do not have the good conservative structure we must benefit from the sharp decay rates of $(\nabla\bfphi_t,\nabla^2 \bfphi)$ to close the estimates. Yet when doing so one cannot rely on an $L^2\to L^4$ bound that would involve a $(1+\tau)^{-1}$ factor. Instead we leave the purely $L^4$-based analysis and estimate the quadratic $\cN_c[\bV,\bfphi]$ in $L^q$ for some $q<2$. 

Let us stress that the current issue differs dramatically from the situation in \cite{MR} where those commutator terms turn out to be subcritical thus asymptotically irrelevant so that they may be estimated in a non sharp way. A similar issue does occur in \cite{JNRZ-conservation} but since, there, one does not need to cope with severe restriction in localization it does not require much attention. To be more explicit on the latter we recall that on one hand $\nabla \bfphi$ does not belong to $L^2$ and that on the other hand in the present multidimensional case involving $L^p$ estimates of $(\bV,\nabla\bfphi)$, or of any of their derivatives, with $p<2$ could result in supercritical decay rates and thus must be avoided in any case. Incidentally we point out that the difficulty to overcome so as to analyze the critical decay dynamics for small solutions of the two-dimensional barotropic compressible Navier-Stokes system is precisely of this kind; see the discussion in \cite{Rodrigues-these}. 

Here the fix is that in the closing of $W^{2,4}$ bounds we need to also involve a sharp control on $(\nabla\bfphi_t,\nabla^2 \bfphi)$ in $L^p$ for some $2\leq p<4$ enabling a control of $\cN_c[\bV,\bfphi]$ in $L^q$ with $1/q=1/4+1/p$ leading to a bound 
\[
\int_0^t\frac{1}{(1+t-\tau)^{\frac1p}}\frac{1}{(1+\tau)^{\frac54-\frac1p}}\dD\tau 
\lesssim\frac{1}{(1+t)^{\frac14}}\,.
\]
With this mind we chose $p=2$ in the above discussion and introduce for any $T>0$ such that $(\bV,\bfphi)$ is defined on $[0,T]$,
\begin{align*}
\zeta(T)&:=
\max_{0\leq t\leq T}(1+t)^{\frac14}\,\|(\bV,\nabla\bfphi,\bfphi_t)(t,\cdot)\|_{W^{2,4}}
+\max_{0\leq t\leq T}(1+t)^{\frac12}\,\|(\cJ\bV,\nabla^2\bfphi,\nabla\bfphi_t)(t,\cdot)\|_{L^2}\,.
\end{align*}
Before moving on, we would like to comment further on the introduction of $\cJ\bV$. Let us first observe that from \eqref{bfphi} it follows that when $t\geq 1$, $\bfphi(t,\cdot)$ is low-frequency so that $\|(\nabla^2\bfphi,\nabla\bfphi_t)(t,\cdot)\|_{L^2}=\|(\cJ\nabla\bfphi,\cJ\bfphi_t)(t,\cdot)\|_{L^2}$. Moreover, when estimating the contribution to $(\cJ\bfphi_t,\cJ\nabla\bfphi)(t,\cdot)$ of conservative terms in integrals over $[t/2,t]$, one cannot leave all the decay on the semigroup part since this would violate the $\alpha<1$ constraint. Instead we distribute the $\cJ$ operator on $(\cN_f[\bV,\bfphi],\cN_s[\bV])$ which is only useful if further decay on $\cJ\bV$ is proved simultaneously.

To begin with, note that a $W^{2,4}$ energy estimate based on \eqref{veq0} and Lemma~\ref{l:energy-estimate} yields for some positive constant $\omega$ and $C$, that if $\zeta(T)$ is defined and sufficiently small, for any $0\leq t\leq T$,
\begin{align*}
\|\bV(t,\cdot)\|_{W^{2,4}}^4
&\leq C\,E_0^4\,\eD^{-\omega\,t}\\
&+C\,
\left(\max_{0\leq t\leq T}(1+t)^{\frac14}\,\left(\|\bV(t,\cdot)\|_{L^4}+\|(\nabla\bfphi,\bfphi_t)(t,\cdot)\|_{W^{2,4}}\right)\right)^4\,\int_0^t \eD^{-\omega\,(t-\tau)}\,\frac{\dD \tau}{(1+\tau)}\,.
\end{align*}
This implies for some positive constant $C$, that if $\zeta(T)$ is defined and sufficiently small
\begin{align*}
\max_{0\leq t\leq T}(1+t)^{\frac14}\,\|\bV(t,\cdot)\|_{W^{2,4}}
\leq C\,E_0+C\,
\max_{0\leq t\leq T}(1+t)^{\frac14}\,\left(\|\bV(t,\cdot)\|_{L^4}+\|(\nabla\bfphi,\bfphi_t)(t,\cdot)\|_{W^{2,4}}\right)\,.
\end{align*}
 
Then proceeding as sketched above to estimate \eqref{bfphi}-\eqref{V} with linear estimates of the previous section one derives for some positive constant $C$, that if $\zeta(T)$ is defined and sufficiently small, for any $0\leq t\leq T$,
\begin{align*}
\|\bV(t,\cdot)\|_{L^4}&+\|(\nabla\bfphi,\bfphi_t)(t,\cdot)\|_{W^{2,4}}\\
&\leq \frac{C\,E_0}{(1+t)^{\frac14}}
+C\,\zeta(T)^2\,\left(
\int_0^t\frac{1}{(1+t-\tau)^{\frac34}}\frac{1}{(1+\tau)^{\frac12}}\dD\tau
+\int_0^t\frac{1}{(1+t-\tau)^{\frac12}}\frac{1}{(1+\tau)^{\frac34}}\dD\tau
\right)\,.
\end{align*}
This implies for some positive constant $C$, that if $\zeta(T)$ is defined and sufficiently small, 
\[
\max_{0\leq t\leq T}(1+t)^{\frac14}\,\left(\|\bV(t,\cdot)\|_{L^4}+\|(\nabla\bfphi,\bfphi_t)(t,\cdot)\|_{W^{2,4}}\right)
\leq C\,E_0+C\,\zeta(T)^2\,.
\]
Inserting this into the above $W^{2,4}$ bound yields for some positive constant $C$, that if $\zeta(T)$ is defined and sufficiently small, 
\[
\max_{0\leq t\leq T}(1+t)^{\frac14}\,\|(\bV,\nabla\bfphi,\bfphi_t)(t,\cdot)\|_{W^{2,4}}
\leq C\,E_0+C\,\zeta(T)^2\,.
\]

It remains to bound $(\cJ\bV,\nabla^2\bfphi,\nabla\bfphi_t)$. This is both by far the most technical part and the part that differs the most from arguments of \cite{JNRZ-conservation,MR}. Let us begin by the elementary remark that both the time-layer contributions arising from $\tchi$ and the $s(0)\cN[\bV(t,\cdot),\bfphi(t,\cdot)]$ contributions to $\bfphi_t(t,\cdot)$ are readily estimated.

At this stage it is convenient to split $\bV$ as $\bV=\bV_{LF}+\bV_{HF}$ according to
\[
\bV_{LF}(t,\cdot):=S_{2}(t) \left[ \bV_{0}+(\bfphi_{0} \cdot \nabla) \ubU \right] + \int_{0}^{t} S_{2}(t-\tau)\cN[\bV(\tau,\cdot),\bfphi(\tau,\cdot)]\dD\tau\,.
\]
Then benefiting from the exponential decay of $S_1$ one derives for some positive constant $C$, that if $\zeta(T)$ is defined and sufficiently small,
\begin{align*}
\max_{0\leq t\leq T}(1+t)^{\frac12}\|\bV_{HF}(t,\cdot)\|_{H^1}
\leq C\,E_0+C\,\zeta(T)^2\,.
\end{align*}
This contains the needed bound on $\|\cJ\bV_{HF}\|_{L^2}$. It also implies through $L^{\frac43}\to L^2$ bounds on $S_2$ and $\d_t s$ that the contribution of $\cN[\bV,\bfphi]-\cN[\bV_{LF},\bfphi]$ to $\|(\cJ\bV_{LF},\nabla^2\bfphi,\nabla\bfphi_t)(t,\cdot)\|_{L^2}$ may be estimated when $0\leq t\leq T$ by a multiple of
\[
(E_0+\zeta(T)^2)^2\,
\int_0^t\frac{1}{(1+t-\tau)^{\frac34}}\frac{1}{(1+\tau)^{\frac34}}\dD\tau 
\lesssim \frac{(E_0+\zeta(T)^2)^2}{(1+t)^{\frac12}}\,.
\]
To apply Lemma~\ref{l:Leibniz} to the bilinear part of $\cN[\bV_{LF},\bfphi]$, let us observe that the contribution of the higher-order --- thus at least cubic --- part of $\cN[\bV_{LF},\bfphi]$ may be estimated similarly through $L^{\frac43}\to L^2$ bounds leading to bounds by a multiple of $\zeta(T)^3\,(1+t)^{-\frac12}$. Likewise $L^{\frac43}\to L^2$ linear bounds are also sufficient to deal with (the rest of) the contribution of $\cN_c[\bV_{LF},\bfphi]$ and yield a bound by a multiple of $\zeta(T)^2\,(1+t)^{-\frac12}$. 

There remains to estimate the contributions of the bilinear parts of $\cN_f[\bV_{LF},\bfphi]$ and $\cN_s[\bV_{LF},\bfphi]$ to $\|(\cJ\bV,\nabla^2\bfphi,\nabla\bfphi_t)\|_{L^2}$. Since, when doing so, we may involve $\|\cJ\bV\|_{L^2}$ but not $\|\cJ\nabla\bV\|_{L^2}$, let us observe that in (the bilinear part of) $\cN_f[\bV_{LF},\bfphi]$ we may factorize all derivatives on $\bV$ up to commutator terms that may be treated as $\cN_c[\bV_{LF},\bfphi]$ was and thus use regularizing effects of $S_2$ and $s$. For the remainders contributions we estimate the $[0,t/2]$-part of the integral through $L^2\to L^2$ linear bounds resulting in a multiple of
\[
\zeta(T)^2\,
\int_0^{t/2}\frac{1}{(1+t-\tau)}\frac{1}{(1+\tau)^{\frac12}}\dD\tau 
\lesssim \frac{\zeta(T)^2}{(1+t)^{\frac12}}\,,
\]
whereas in the $[t/2,t]$-part of the integral we distribute the $\cJ$ operator on the bilinear low-Floquet terms along Lemma~\ref{l:Leibniz} with $m=2$ and then use $L^{\frac32}\to L^2$ linear bounds deriving an estimate by a multiple of $\zeta(T)^2\,(1+t)^{-\frac12}$.
 
This completes the proof that for some positive constant $C$, if $\zeta(T)$ is defined and sufficiently small,
\[
\zeta(T)\leq C\,(E_0+\zeta(T)^2)\,(1+E_0+\zeta(T)^2)\,.
\]
In turn this implies for the same $C$, that provided that $E_0$ is sufficiently small, if the solution exists on $[0,T]$ and satisfies $\zeta(T)\leq 3C\,E_0$ then it also satisfies $\zeta(T)\leq 2C\,E_0$. From a continuity argument one deduces that the solution is defined on $\R_+$ and satisfies 
\[
\sup_{t\geq 0}\,(1+t)^{\frac14}\,\|(\bV,\nabla\bfphi,\bfphi_t)(t,\cdot)\|_{W^{2,4}}
+\sup_{t\geq 0}\,(1+t)^{\frac12}\,\|(\cJ\bV,\nabla^2\bfphi,\nabla\bfphi_t)(t,\cdot)\|_{L^2}
\,\leq\, 2C\,E_0\,.
\]

To complete the proof of Theorem~\ref{theo:critmain} there only remains to notice that this is sufficient to derive the claimed $L^p$ bounds. Incidentally we point out that we have decided to keep the derivation of those $L^p$ bounds simple but a finer analysis in the spirit of the above $\cJ$-analysis would remove the $\ln(2+t)$ factor of the $L^\infty$ bound.

\section{Maximal decay regime}

We now focus on the proof of Theorem~\ref{theo:submain}. For the sake of comparison we point out that the closest statement in \cite{MR} is \cite[Theorem~1.2]{MR}, which is a subcritical but not a maximal decay regime there.

\subsection*{Linear estimates}

Most of the elements of the proof of Theorem~\ref{theo:critmain} may be reused as they are but we need to revisit the linear estimates for initial phase modulations provided by Proposition~\ref{prop:linphasedecay} and Lemma~\ref{l:boundarylayer}. Phases are here obtained from their gradient through $\bfphi=-(-\Delta)^{-1}\,\Div(\nabla\bfphi)$. Slight variations on the original proofs provide the relevant replacements. 

\bpr
Assume \cond1-\cond2-\cond3.
\begin{enumerate}
\item There exists $\theta'>0$, such that, for any $(s,s')\in\R_+$ such that $s'\leq s$ and any $(p_0,p_1)$ such that $2<p_0<p_1<\infty$, there exists $C_{p_0,p_1,s,s'}$ such that for any $t>0$, and any $p_0\leq p\leq p_1$,
\begin{align*}
\|S_{1}(t)[(\bfphi \cdot \nabla) \ubU]\|_{W^{s,p}}
&\leq \frac{C_{p_0,p_1,s,s'}}{(\min(\{1,t\}))^{\frac{(s-(s'+1))}{2}}}\,\eD^{-\theta'\,t}
\,\|\nabla\bfphi\|_{L^1\cap H^{s'}}\,.
\end{align*}
\item For any $\alpha\in\N^2$ and any $s\in\R_+$, there exists $C_{\alpha,s}$ such that for any $2 \leq p \leq +\infty$ and any $t\geq0$
\begin{align*}
\|\cJ^\alpha S_{2}(t)[(\bfphi \cdot \nabla) \ubU]\|_{W^{s,p}}
&\leq \frac{C_{\alpha,s}}{(1+t)^{\frac{|\alpha|}{2}+1-\frac{1}{p}}}\,\|\nabla\bfphi\|_{L^1}\,.
\end{align*}
\item For any $\alpha\in\N^2$, any $\beta\in\N^2$, any $\ell\in\N$ and any $2 \leq p \leq +\infty$ such that $|\alpha|+|\beta|+\ell\geq 1$ or $p>2$, there exists $C_{p,\alpha,\beta,\ell}$ such that for any $t\geq0$
\begin{align*}
\|\,\d_\bfx^\alpha\,\d_t^\ell\,\cJ^\beta\,s(t)[(\bfphi \cdot \nabla) \ubU]\|_{L^p}
&\leq \frac{C_{p,\alpha,\ell}}{(1+t)^{\frac{|\alpha|+|\beta|+\ell}{2}+\frac12-\frac{1}{p}}}\,\|\nabla\bfphi\|_{L^1}\,.
\end{align*}
\end{enumerate}
\epr

For the sake of comparison we have kept the $\cJ$ operator in the above statement but it is of no use in the proof of Theorem~\ref{theo:submain}.  

\bl
Assume \cond1-\cond2-\cond3. For any $\alpha\in\N^2$, any $\ell\in\N$ and any $2 \leq p \leq +\infty$ such that $p<\infty$ if $|\alpha|\geq 1$, there exists $C_{p,\alpha}$ such that for any $t\geq 0$
\begin{align*}
\|\,\d_\bfx^\alpha\,\left(s(t)[(\bfphi \cdot \nabla)\ubU]-\bfphi\right)\|_{L^p}
&\leq C_{p,\alpha}\,\|\nabla\bfphi\|_{L^1\cap W^{(|\alpha|-1)_+,p}}
\,\begin{cases}
\,(1+t)^{\frac12\,\left(\frac2p-|\alpha|\right)_+}&\quad\textrm{if } |\alpha|-\tfrac2p\neq 0\\
\ln(2+t)&\quad\textrm{otherwise }
\end{cases}\,.
\end{align*}
\el

\subsection*{Nonlinear high-frequency damping estimate}

To propagate the better localization we complete the $L^4$-dissipation estimates of Lemma~\ref{l:energy-estimate} with standard $L^2$-based energy estimates that we omit to state. Actually it would be possible here to obtain a statement similar to Theorem~\ref{theo:submain} with purely $L^2$-based damping estimates.

\subsection*{Proof of Theorem~\ref{theo:submain}}

As in the proof of Theorem~\ref{theo:critmain} we first close a first round of bounds by a continuity argument and then prove the remaining bounds using the first ones. A convenient set of norms in order to achieve the nonlinear closing part is 
\begin{align*}
\zeta(T)&:=
\max_{0\leq t\leq T}(1+t)^{\frac34}\,\|(\bV,\nabla\bfphi,\bfphi_t)(t,\cdot)\|_{W^{2,4}}
+\max_{0\leq t\leq T}(1+t)^{\frac54}\,\|(\nabla^2\bfphi,\nabla\bfphi_t)(t,\cdot)\|_{L^4}\\
&\quad
+\max_{0\leq t\leq T}(1+t)^{\frac12}\,\|(\bV,\nabla\bfphi,\bfphi_t)(t,\cdot)\|_{H^{2}}
+\max_{0\leq t\leq T}(1+t)\,\|(\nabla^2\bfphi,\nabla\bfphi_t)(t,\cdot)\|_{L^2}
\,.
\end{align*}
The mixed $L^2/L^4$ framework is designed to be able to place nonlinear terms both in $L^1$ and in $L^2$. Note that we do not need to involve the $\cJ$ operator mostly because there is less difficulty to propagate even high decay rates in a subcritical regime.

Before giving some more concrete elements we point out that in subcritical regimes one make an intensive use of 
\begin{align*}
\int_0^{t/2}\frac{1}{(1+t-\tau)^\alpha}\frac{1}{(1+\tau)^\beta}\dD\tau &\lesssim\frac{1}{(1+t)^{\beta}}\,,& 0\leq \alpha\,,\quad 1<\beta\,,\\
\int_{t/2}^t\frac{1}{(1+t-\tau)^\alpha}\frac{1}{(1+\tau)^\beta}\dD\tau &\lesssim\frac{1}{(1+t)^{\beta}}\,,& 1<\alpha\,,\quad 0\leq \beta\,,
\end{align*}
and
\begin{align*}
\int_0^{t/2}\frac{1}{(1+t-\tau)^\alpha}\frac{1}{(1+\tau)}\dD\tau &\lesssim\frac{\ln(2+t)}{(1+t)^{\alpha}}\,,& 0\leq \alpha\,,\\
\int_{t/2}^t\frac{1}{(1+t-\tau)}\frac{1}{(1+\tau)^\beta}\dD\tau &\lesssim\frac{\ln(2+t)}{(1+t)^{\beta}}\,,& 0\leq \beta\,.
\end{align*}

To evaluate the $[0,t/2]$-part of the integrals we use $L^1\to L^p$ bounds, with $p=2,4$, resulting in bounds by multiples of the following quantities 
\begin{enumerate}
\item when estimating the contribution of flux and source terms to $\|(\bV,\nabla\bfphi,\bfphi_t)(t,\cdot)\|_{L^p}$
\[
\zeta(T)^2\,\int_0^{t/2}\frac{1}{(1+t-\tau)^{\frac12+1-\frac1p}}\frac{1}{(1+\tau)}\dD\tau 
\lesssim\frac{\zeta(T)^2\,\ln(2+t)}{(1+t)^{\frac12+1-\frac1p}}
\lesssim\frac{\zeta(T)^2}{(1+t)^{1-\frac1p}}
\]
\item when estimating the contribution of commutator terms to $\|(\bV,\nabla\bfphi,\bfphi_t)(t,\cdot)\|_{L^p}$
\[
\zeta(T)^2\,\int_0^{t/2}\frac{1}{(1+t-\tau)^{1-\frac1p}}\frac{1}{(1+\tau)^{\frac32}}\dD\tau 
\lesssim\frac{\zeta(T)^2}{(1+t)^{1-\frac1p}}
\]
\end{enumerate}
and similarly with extra $(1+t-\tau)^{-\frac12}$ and $(1+t)^{-\frac12}$ factors for contributions to $\|(\nabla^2\bfphi,\nabla\bfphi_t)(t,\cdot)\|_{L^p}$. 

In turn, to evaluate the $[t/2,t]$-part of the integrals we use $L^2\to L^p$ bounds, with $p=2,4$, resulting in bounds by multiples of the following quantities 
\begin{enumerate}
\item when estimating the contribution of flux and source terms to $\|(\nabla^2\bfphi,\nabla\bfphi_t)(t,\cdot)\|_{L^p}$
\[
\zeta(T)^2\,\int_{t/2}^t\frac{1}{(1+t-\tau)^{\frac12+1-\frac1p}}\frac{1}{(1+\tau)^{\frac32}}\dD\tau 
\lesssim\frac{\zeta(T)^2}{(1+t)^{\frac32}}
\lesssim\frac{\zeta(T)^2}{(1+t)^{\frac32-\frac1p}}
\]
\item when estimating the contribution of commutator terms to $\|(\nabla^2\bfphi,\nabla\bfphi_t)(t,\cdot)\|_{L^p}$
\[
\zeta(T)^2\,\int_{t/2}^t\frac{1}{(1+t-\tau)^{1-\frac1p}}\frac{1}{(1+\tau)^2}\dD\tau 
\lesssim\frac{\zeta(T)^2}{(1+t)^{\frac32-\frac1p}}
\]
\end{enumerate}
and similarly with extra $(1+t-\tau)^{\frac12}$ and $(1+t)^{\frac12}$ factors for contributions to $\|(\bV,\nabla\bfphi,\bfphi_t)(t,\cdot)\|_{L^p}$.

This is essentially sufficient to close the continuity argument and then prove the rest of the bounds. Again we have decided not to put too much energy in removing spurious $\ln(2+t)$ factors but we'd like to point out that one way to discard those here is to prove and use estimates that are deteriorated by dispersion, the subcriticality compensating for the dispersive loss. We refer to \cite{MR} for related analyses.

In subcritical regimes one expects nonlinear contributions to be asymptotically irrelevant. Yet our proof does not establish this claim in the present situation. Further work would be needed to prove it and again we refer to \cite{MR} for related analyses.

\appendix
\section{Extensions to multiD}\label{s:3D}

We now turn our attention to higher dimensional cases. When going to higher dimensions the variety of scenarios to consider increases both in types of periodic waves and in range of possible localizations between maximal and critical decay. We focus here on genuinely $d$-dimensional periodic waves in dimension $d$ and critically localized perturbations. We expect the treatment of more localized perturbations to follow by an easier version of the argument and we refer the reader to \cite{JZ-multiD} for an example of analysis near a $1$D wave in dimension $d\geq3$ under subcritical perturbations.

The translation of the problem in the multidimensional framework is immediate. In the set of assumptions the only significant change is that in \cond{1}, $0$ is assumed to be an eigenvalue of $L_{\b0}$ of multiplicity $d+r$. From the point of view of time decay, any part of the decay due to cancellation at low-Floquet exponents, arising for instance by taking derivatives on phases, is insensitive to dimension whereas variations in Lebesgue indices do feel the dimension, $1/q$ and $1/p$ being now preceded by a $d/2$ slope. We choose to encode the critical localization assumption as $(-\Delta)^{\frac{d}{2}}\bfphi_0\in L^1 (\R^d)$. Note that $\nabla^d\bfphi_0\in L^1 (\R^d)$ or, when $d$ is odd, $(-\Delta)^{\frac{d-1}{2}}\nabla\bfphi_0\in L^1 (\R^d)$ would also be reasonable nonequivalent choices. 

\subsection*{3D}

Despite immediate similarities in arbitrary dimension the only dimension where our analysis transfers almost word by word is dimension $d=3$, resulting in the following result. 

\bt[d=3]
Set $d=3$ and consider a $d$-D periodic wave in $\R^d$ satisfying the $d$-D version of diffuse spectral stability \cond1-\cond2-\cond3. There exist $\eps_0>0$ and $C>0$ such that if for some sublinear\footnote{By this, we mean that $\bfphi_0$ may differ from $(-\Delta)^{-\frac{d}{2}}((-\Delta)^{\frac{d}{2}} \bfphi_0)$ only by a constant function.} $\bfphi_0$
\[
E_0:=\|\bW_{0}(\cdot-\bfphi_{0})-\ubU \|_{(H^{\min(\{2;d\})} \cap W^{2,4})(\R^d;\R^n)}+\|(-\Delta)^{\frac{d}{2}}\bfphi_{0}\|_{(L^{1} \cap W^{3-d,4})(\R^d;\R^d)}\,\leq\,\eps_0
\]
then, there exist a unique global solution with initial datum $\bW_{0}$ and a phase shift $\bfphi$ with $\bfphi(0,\cdot)= \bfphi_{0}$ such that, for any $t\geq0$,
\begin{align*}
\|\bW(t,\cdot-\bfphi(t,\cdot))-\ubU\|_{W^{2,4}(\R^d;\R^n)}
+\|\nabla\bfphi(t,\cdot)\|_{W^{2,4}(\R^d;\cM_d(\R))}
+\|\d_t\bfphi(t,\cdot)\|_{W^{2,4}(\R^d;\R^d)}
&\leq \frac{C\,E_{0}}{(1+t)^{\frac12-\frac{d}{2}\frac{1}{4}}}\,.
\end{align*}
Furthermore, with $p_d$ defined by $1/p_d:=1/2-1/d$, constants independent of $(\bW_0,\bfphi_0)$ and no further restriction on $E_0$, 
\begin{enumerate}
\item for any $t\geq0$, 
\begin{align*}
\|\bW(t,\cdot-\bfphi(t,\cdot))-\ubU\|_{L^{p_q}(\R^d;\R^n)}
+\|\nabla\bfphi(t,\cdot)\|_{L^\infty(\R^d;\cM_d(\R))}
+\|\d_t\bfphi(t,\cdot)\|_{L^\infty(\R^d;\R^d)}
&\leq C\,E_0\,\frac{\ln(2+t)}{(1+t)^{\frac{1}{2}-\frac{d}{2}\frac{1}{p_d}}}\,;
\end{align*}
\item for any $d<p_{0}<q_{0}<p_d$, there exists a constant $C_{p_{0},q_{0}}>0$, such that for any $p \in [p_{0},q_{0}]$, and any $t\geq0$
\[
\|\bW(t,\cdot-\bfphi(t,\cdot))-\ubU\|_{L^p(\R^d;\R^n)}
+\|\nabla\bfphi(t,\cdot)\|_{L^p(\R^d;\cM_d(\R))}
+\|\d_t\bfphi(t,\cdot)\|_{L^p(\R^d;\R^d)}
\,\leq \frac{C_{p_{0},q_{0}}\,E_{0}}{(1+t)^{\frac{1}{2}-\frac{d}{2}\frac{1}{p}}}\,.
\]
\end{enumerate}
\et

Despite the fact that the previous statement only holds in dimension $d=3$ we have placed the general letter $d$ in most places to facilitate the later comparison with the general case. In $L^p$ bounds, we stress that $p>d$ is a hard constraint, the functions of interest do not belong to $L^d$, whereas $p<p_d$ is a soft constraint, our method of proof does not give sharp decay rates when $p>p_d$ because of the saturation of some time-integration constraints.

Since the proof of the foregoing theorem is essentially identical to the one of Theorem~\ref{theo:critmain} we only discuss parts where in principle dimension could matter. 
\begin{enumerate}
\item To begin with we note that, in dimension $3$, the two Sobolev embeddings used in the proof still hold, namely $H^1$ embeds in $L^4$ and $W^{2,4}$ embeds in $W^{1,\infty}$.
\item More importantly, in dimension $3$, $(-\Delta)^{\frac{3}{2}}\bfphi_{0}\in L^1\cap L^4$ does imply both $\nabla \bfphi\in W^{2,4}$ and $\nabla^2\bfphi\in L^2$. This follows from $\bfxi\mapsto \|\bfxi\|^{-(d-1)}\in L^{4'}$ and $\bfxi\mapsto \|\bfxi\|^{-(d-2)}\in L^{2'}$.
\item All the time-integration constraints are satisfied. With this respect we point out that due to the monotonicity in $d$ of the involved time powers only the powers of $(1+t-\tau)$ in the $[t/2,t]$ parts of the integrals must be checked. To mention only the worst constraints, with $d=3$ one still has
\begin{align*}
\frac{d}{2}\left(\frac12-\frac14\right)+\frac12&<1\,,&
\frac{d}{2}\left(\frac34-\frac14\right)&<1\,.
\end{align*}
\item The $L^2$-based estimates used to bound contributions associated with $S_1$ are better than the aimed $L^4$ estimates. Indeed with $d=3$, $1-\frac{d}{2}\frac12\geq\frac12-\frac{d}{2}\frac14$.
\end{enumerate}

\subsection*{General dimension}

Let us now return to the discussion of the general multiD case and examine whether one may modify some of the arbitrary choices made in the proof of Theorem~\ref{theo:critmain} so as to relax the constraint $d\leq 3$. We need to pick $r_1$ and $r_2$ to place $(\bV,\nabla \bfphi)\in W^{s,r_1}$ (for some sufficiently large $s$) and $(\cJ\bV,\nabla^2 \bfphi)\in L^{r_2}$. The first hard constraint is that this should be compatible with assumptions on $(-\Delta)^{\frac{d}{2}} \bfphi_0$, namely with $\bfxi\mapsto \|\bfxi\|^{-(d-1)}\in L^{r_1'}$ and $\bfxi\mapsto \|\bfxi\|^{-(d-2)}\in L^{r_2'}$. Those are equivalent to $r_1>d$ and $r_2>d/2$. Yet this would place nonlinear terms only in $L^{r_1/2}$ with $r_1/2>2$ when $d\geq4$. This is incompatible with the kind of linear bounds that we have proved here by resorting mostly to Hausdorff-Young inequalities.  

A possible way out would be to distribute powers of $\cJ$, this time not to balance time decay but to enhance localization. Note in particular that for $0\leq \ell<d$, one may derive $\nabla^\ell\bfphi\in L^{r_\ell}$ when $r_\ell>d/\ell$.

Since the constraint $d\leq 3$ already covers most of cases of practical interest we restrain ourselves from delving into this avenue and leave this for possible future work. In any case we recall that we expect that with the techniques of the present contribution the maximal decay regime could be covered in any dimension in a more elementary way.


\begin{thebibliography}{JNRZ13b}

\bibitem[AR22]{Audiard-Rodrigues}
C.~Audiard and L.~M. Rodrigues.
\newblock About plane periodic waves of the nonlinear {S}chr\"odinger
  equations.
\newblock {\em Bull. Soc. Math. France}, 150(1):111--207, 2022.

\bibitem[BGNR14]{BGNR}
S.~Benzoni-Gavage, P.~Noble, and L.~M. Rodrigues.
\newblock Slow modulations of periodic waves in {H}amiltonian {PDE}s, with
  application to capillary fluids.
\newblock {\em J. Nonlinear Sci.}, 24(4):711--768, 2014.

\bibitem[BGS07]{Benzoni-Serre}
S.~Benzoni-Gavage and D.~Serre.
\newblock {\em Multidimensional hyperbolic partial differential equations}.
\newblock Oxford Mathematical Monographs. The Clarendon Press, Oxford
  University Press, Oxford, 2007.
\newblock First-order systems and applications.

\bibitem[dR24]{dR}
B.~de~Rijk.
\newblock Nonlinear stability and asymptotic behavior of periodic wave trains
  in reaction-diffusion systems against {$C_{\rm ub}$}-perturbations.
\newblock {\em Arch. Ration. Mech. Anal.}, 248(3):Paper No. 36, 53, 2024.

\bibitem[GR25]{GR}
L.~Gar\'enaux and L.~M. Rodrigues.
\newblock Convective stability in scalar balance laws.
\newblock {\em Differential Integral Equations}, 38(1-2):71--110, 2025.

\bibitem[HZ95]{Hoff_Zumbrun-NS_compressible_pres_de_zero}
D.~Hoff and K.~Zumbrun.
\newblock Multi-dimensional diffusion waves for the {N}avier-{S}tokes equations
  of compressible flow.
\newblock {\em Indiana Univ. Math. J.}, 44(2):603--676, 1995.

\bibitem[JNRZ13a]{JNRZ-RD1}
M.~A. Johnson, P.~Noble, L.~M. Rodrigues, and K.~Zumbrun.
\newblock Nonlocalized modulation of periodic reaction diffusion waves:
  nonlinear stability.
\newblock {\em Arch. Ration. Mech. Anal.}, 207(2):693--715, 2013.

\bibitem[JNRZ13b]{JNRZ-RD2}
M.~A. Johnson, P.~Noble, L.~M. Rodrigues, and K.~Zumbrun.
\newblock Nonlocalized modulation of periodic reaction diffusion waves: the
  {W}hitham equation.
\newblock {\em Arch. Ration. Mech. Anal.}, 207(2):669--692, 2013.

\bibitem[JNRZ14]{JNRZ-conservation}
M.~A. Johnson, P.~Noble, L.~M. Rodrigues, and K.~Zumbrun.
\newblock Behavior of periodic solutions of viscous conservation laws under
  localized and nonlocalized perturbations.
\newblock {\em Invent. Math.}, 197(1):115--213, 2014.

\bibitem[JZ10]{JZ-multiD}
M.~A. Johnson and K.~Zumbrun.
\newblock Nonlinear stability of periodic traveling wave solutions of systems
  of viscous conservation laws in the generic case.
\newblock {\em J. Differential Equations}, 249(5):1213--1240, 2010.

\bibitem[JZ11]{JZ}
M.~A. Johnson and K.~Zumbrun.
\newblock Nonlinear stability of spatially-periodic traveling-wave solutions of
  systems of reaction-diffusion equations.
\newblock {\em Ann. Inst. H. Poincar\'e{} C Anal. Non Lin\'eaire},
  28(4):471--483, 2011.

\bibitem[KR16]{KR}
B.~Kabil and L.~M. Rodrigues.
\newblock Spectral validation of the {W}hitham equations for periodic waves of
  lattice dynamical systems.
\newblock {\em J. Differential Equations}, 260(3):2994--3028, 2016.

\bibitem[MR24]{MR}
B.~Melinand and L.~M. Rodrigues.
\newblock Phase sinks and sources around two-dimensional periodic-wave
  solutions of reaction-diffusion-advection systems.
\newblock {\em arXiv preprint arXiv:2408.14869}, 2024.

\bibitem[NR13]{Noble-Rodrigues}
P.~Noble and L.~M. Rodrigues.
\newblock Whitham's modulation equations and stability of periodic wave
  solutions of the {K}orteweg-de {V}ries-{K}uramoto-{S}ivashinsky equation.
\newblock {\em Indiana Univ. Math. J.}, 62(3):753--783, 2013.

\bibitem[{\SortNoop{Rijk}}RR]{BdR-R}
B.~{\SortNoop{Rijk}}de~Rijk and L.~M. Rodrigues.
\newblock Secondary instabilities in reaction-diffusion-advection systems:
  far-from-equilibrium pattern formation due to transverse destabilization.
\newblock Forthcoming.

\bibitem[Rod07]{Rodrigues-these}
L.~M. Rodrigues.
\newblock {\em Comportement en temps long des fluids visqueux bidimensionnels}.
\newblock PhD thesis, Universit\'e Grenoble 1, 2007.
\newblock In French.

\bibitem[Rod09]{Rodrigues-compressible}
L.~M. Rodrigues.
\newblock Vortex-like finite-energy asymptotic profiles for isentropic
  compressible flows.
\newblock {\em Indiana Univ. Math. J.}, 58(4):1747--1776, 2009.

\bibitem[Rod13]{R}
L.~M. Rodrigues.
\newblock {\em Asymptotic stability and modulation of periodic wavetrains,
  general theory \& applications to thin film flows}.
\newblock Habilitation {\`a} diriger des recherches, Universit\'e Lyon 1, 2013.

\bibitem[Rod15]{R_Roscoff}
L.~M. Rodrigues.
\newblock Space-modulated stability and averaged dynamics.
\newblock {\em Journ\'ees \'Equations aux d\'eriv\'ees partielles},
  2015(8):1--15, 2015.

\bibitem[Rod18]{R-linKdV}
L.~M. Rodrigues.
\newblock Linear asymptotic stability and modulation behavior near periodic
  waves of the {K}orteweg--de {V}ries equation.
\newblock {\em J. Funct. Anal.}, 274(9):2553--2605, 2018.

\bibitem[SSSU12]{SSSU}
B.~Sandstede, A.~Scheel, G.~Schneider, and H.~Uecker.
\newblock Diffusive mixing of periodic wave trains in reaction-diffusion
  systems.
\newblock {\em J. Differential Equations}, 252(5):3541--3574, 2012.

\bibitem[Sat79]{Sattinger_bifurcation}
D.~H. Sattinger.
\newblock {\em Group-theoretic methods in bifurcation theory}, volume 762 of
  {\em Lecture Notes in Mathematics}.
\newblock Springer, Berlin, 1979.
\newblock With an appendix entitled ``How to find the symmetry group of a
  differential equation'' by Peter Olver.

\bibitem[Sch96]{Schneider}
G.~Schneider.
\newblock Diffusive stability of spatial periodic solutions of the
  {S}wift-{H}ohenberg equation.
\newblock {\em Comm. Math. Phys.}, 178(3):679--702, 1996.

\end{thebibliography}

\end{document}